\documentclass[12pt]{article}
\usepackage{amsfonts,amssymb,amsmath,amsthm}

\usepackage{xcolor}

\usepackage{graphicx}

\newcommand{\R}{{\mathbb R}}

\newcommand{\Z}{{\bf Z}}

\newcommand{\re}{{\rm e}}

\newcommand{\Fix}{{\rm Fix}}
\newcommand{\be}{{\bf e}}
\newcommand{\br}{{\bf r}}
\newcommand{\bw}{{\bf w}}
\newcommand{\bc}{{\bf c}}

\newcommand{\bx}{{\bf x}}
\newcommand{\bq}{{\bf q}}
\newcommand{\bh}{{\bf h}}
\newcommand{\bv}{{\bf v}}

\newcommand{\bt}{{\bf t}}
\newcommand{\bz}{{\bf z}}

\newcommand{\by}{{\bf y}}
\newcommand{\bu}{{\bf u}}

\newcommand{\rd}{{\rm d}}

\newcommand{\cF}{{\cal F}}
\newcommand{\cX}{{\cal X}}

\newcommand{\cP}{{\cal P}}

\newcommand{\cW}{{\cal W}}
\newcommand{\cN}{{\cal N}}

\newcommand{\td}{{\mathfrak{d}}}

\newtheorem{theorem}{Theorem}[section]
\newtheorem{lemma}[theorem]{Lemma}
\newtheorem{Cy}[theorem]{Corollary}
\newtheorem{definition}[theorem]{Definition}
\theoremstyle{definition}
\newtheorem{remark}[theorem]{Remark}
\newtheorem{example}{Example}

\begin{document}
\begin{center}
\mbox{}
\bigskip

\noindent{\large \bf Asymptotic stability of robust heteroclinic networks}

\bigskip

{\small \bf Olga Podvigina}\\
{\footnotesize Institute of Earthquake Prediction Theory\\
and Mathematical Geophysics,\\
84/32 Profsoyuznaya St, 117997 Moscow, Russian Federation \\
email: olgap@mitp.ru}

\bigskip
{\small \bf Sofia B.S.D. Castro$^*$}\\
{\footnotesize Centro de Matem\'atica and Faculdade de Economia\\
Universidade do Porto\\
Rua Dr. Roberto Frias, 4200-464 Porto, Portugal\\
email: sdcastro@fep.up.pt}

\bigskip

{\small \bf Isabel S. Labouriau }\\
{\footnotesize Centro de Matem\'atica da Universidade do Porto\\
Rua do Campo Alegre 687, 4169-007 Porto, Portugal\\
email: islabour@fc.up.pt}
 \end{center}

\medskip

\noindent $^*$ corresponding author

\medskip

\noindent {\bf Keywords:} heteroclinic network; asymptotic stability; symmetry

\begin{abstract}
We provide conditions guaranteeing that certain classes of robust heteroclinic networks are
asymptotically stable.

We study the asymptotic stability of ac-networks --- robust heteroclinic networks that exist in  smooth $\Z^n_2$-equivariant dynamical systems defined in the positive orthant
of $\R^n$.
Generators of the group $\Z^n_2$ are
the transformations that change the sign of one of the spatial coordinates. The
ac-network
 is a union of hyperbolic equilibria and connecting trajectories, where all equilibria belong to the coordinate axes (not more than one equilibrium per axis) with unstable manifolds of dimension one or two.
The classification of ac-networks is carried out by describing all possible types of associated graphs.

We prove sufficient conditions for asymptotic stability of
 ac-networks. The proof is given as a series of theorems and lemmas that are applicable to the ac-networks and to more general types of networks. Finally, we apply these results to
discuss the
 asymptotic stability of several examples of heteroclinic networks.
\end{abstract}

\section{Introduction}

A large number of examples of asymptotically stable heteroclinic cycles
can be found in the literature. See, for instance,
Busse and Heikes \cite{bh80},
Jones and Proctor \cite{jp87},
Guckenheimer and Holmes \cite{gh88},
Hofbauer and So \cite{hs94},
Krupa and Melbourne \cite{KrupaMelbourne1} and \cite{KrupaMelbourne2},
Feng \cite{Feng1998},
Postlethwaite \cite{pos2010}
Podvigina \cite{op12} and \cite{op13},
Lohse \cite{Lohse15} or
Podvigina and Chossat \cite{pc16}
Conditions for asymptotic stability for
heteroclinic cycles in particular systems, or  in certain classes in
low-dimensional systems have been known for a long time, see
\cite{bh80,jp87,gh88,hs94,Feng1998,pos2010}
for the former
and  \cite{KrupaMelbourne1,KrupaMelbourne2,op12,op13,pc16}
as well as Garrido-da-Silva and Castro \cite{GdSC2017}
 for the latter, respectively.

With heteroclinic networks the situation is completely different --- there
are just a few instances of networks whose asymptotic stability was proven. See Kirk {\em et al.} \cite{kpr}
or Afraimovich {\em et al.} \cite{afraimovich16}.
One can point out two reasons for the rareness of such heteroclinic networks
in literature. The first is the rather restrictive necessary conditions,
 stating that a
compact asymptotically stable heteroclinic network contains the unstable
manifolds of all its nodes as subsets. In  Podvigina {\em et al.} \cite{pcl18}  this was proven for
networks where all nodes are equilibria.
Applied to networks with
one-dimensional connections these conditions rule out their asymptotic
stability for a large variety of systems, e.g., equivariant
or population dynamics. The second reason  is the complexity
of the problem. Given a heteroclinic cycle, the derivation of conditions
for asymptotic stability involves the construction of a return map around the
cycle which typically is a highly non-trivial problem. Existence of various
walks along a network that can be followed by nearby trajectories makes
the study of the stability of networks much more difficult than that of cycles.

A heteroclinic cycle is a union of nodes and connecting trajectories and
a heteroclinic network is a union of heteroclinic cycles.
Heteroclinic cycles or networks do not exist in a generic dynamical
system, because small perturbations break connections between saddles.
However, they may exist in systems where some constraints are imposed
and be robust with respect to a constrained perturbation due to the
presence of invariant subspaces. Typically, the constraints create flow-invariant subspaces where the connection is of saddle-sink type.
For instance, in a $\Gamma$-equivariant system
the fixed-point subspace of a subgroup of $\Gamma$ is flow-invariant. In population
dynamics modelled by systems in $\R^n_+$, the ``extinction subspaces'', which are the
Cartesian hyperplanes, are flow-invariant. Similarly, such hyperplanes are
invariant subspaces  for systems on a simplex, a usual state space in evolutionary game theory. For coupled cells or
oscillators, existence of flow-invariant subspaces follows from the prescribed
patterns of interaction between the components.

Heteroclinic networks have very different levels of complexity.
A classification of the least complex networks has been proposed by Krupa and Melbourne \cite{KrupaMelbourne1} (simple), Podvigina and Chossat \cite{pc15}, \cite{pc16} (pseudo-simple) and Garrido-da-Silva and Castro \cite{GdSC2017} (quasi-simple).

In this paper we study heteroclinic networks emerging in a smooth
$\Z_2^n$-equivariant dynamical system defined in the positive
orthant $\R^n_+=\{\bx\in\R^n~:~\min_{1\le k\le n}x_k\ge0\}$.
The group $\Z_2^n$ is the group generated by the transformations
that change the sign of one of the spatial coordinates. The networks
that we study, that we call {\em ac-networks}, are
comprised of hyperbolic  equilibria and connecting trajectories, where all equilibria
belong to the coordinate axes (not more than one equilibrium per axis) with
unstable manifold of dimension one or two.

We classify ac-networks by describing the possible structure of the associated graphs.
Graphs are often employed for visualisation and/or study of heteroclinic
networks. The relations between graphs and networks are discussed in
a number of papers (see Ashwin and Postlethwaite \cite{AshPos13} and Field \cite{fie15}), with particular attention being given to the construction of a dynamical system
that has a heteroclinic network with a given graph.
Similarly to \cite{AshPos13}, an ac-network can be realised as a network on a simplex.
Note that networks obtained by the simplex method of \cite{AshPos13} are
not necessarily ac-networks since the method does not require
the unstable manifold of each equilibrium to be entirely contained in the
network. (See Definitions \ref{clean}
and \ref{acn} for precise statements.)

\begin{figure}[h]
\begin{center}
\parbox{0.45\textwidth}
{\begin{center}
\includegraphics[width=0.3\textwidth]{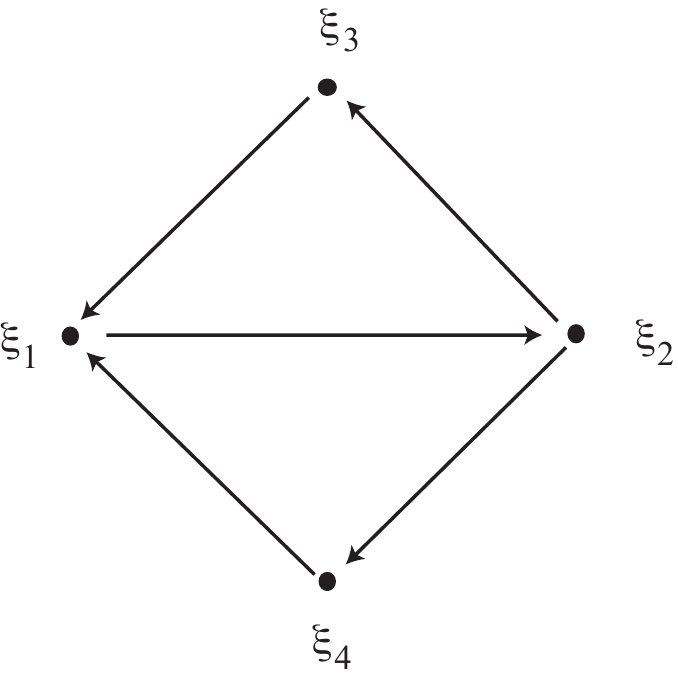}\\
(a)
\end{center}}
\qquad
\parbox{0.45\textwidth}
{\begin{center}
\includegraphics[width=0.3\textwidth]{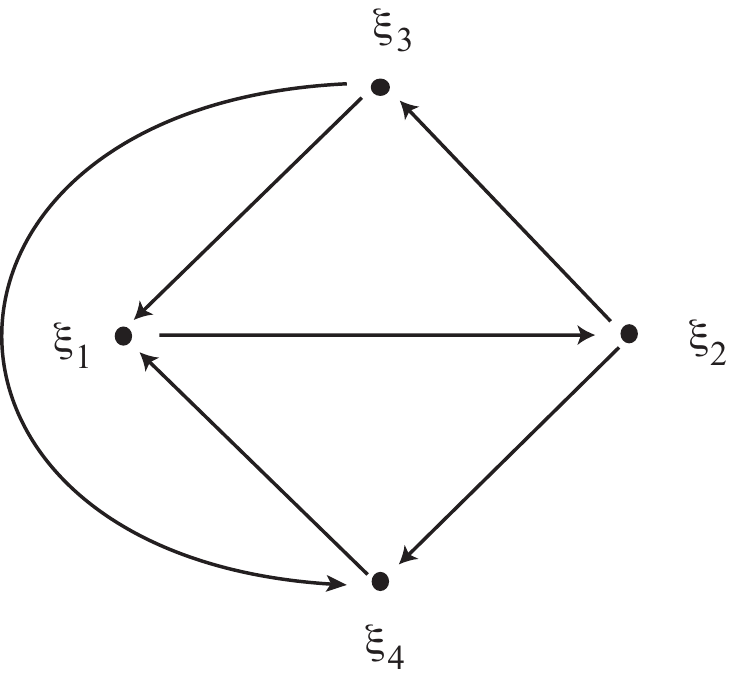}\\
(b)
\end{center}}
\end{center}
\caption{
(a) The Kirk and Silber \cite{KS94} network is not an ac-network, but as shown
in  (b)  in may become an ac-network  by adding
the full sets of connections from $\xi_3$ to $\xi_4$
and from $\xi_2$ to $\xi_4$, as done in \cite{Br,CL16}.
Note that the unstable
manifold of $\xi_2$ is now bounded by the connection from $\xi_3$ to $\xi_4$.
Hence points in the unstable manifold of $\xi_2$ are taken either to $\xi_3$
or to $\xi_4$.
}\label{figKS}
\end{figure}

A very well known network realised by the simplex method is that presented in Kirk and Silber \cite{KS94} and reproduced here in
Figure~\ref{figKS}~(a).
The unstable manifold of $\xi_2$ is not contained in the network.
By  adding
the full sets of connections from $\xi_3$ to $\xi_4$
and from $\xi_2$ to $\xi_4$
to the network we obtain an ac-network (see Figure~\ref{figKS}~(b)).
The stability of this network was discussed in Brannath \cite{Br}
and Castro and Lohse \cite{CL16}.

We prove sufficient conditions for the asymptotic stability of heteroclinic
networks. As usual, we approximate the behaviour of trajectories by local
and global maps. The derived estimates for local maps near $\xi_j$ involve the
exponents $\rho_j$. For compact networks global maps are linear with
constants that are bounded from above. We show that a network is
asymptotically stable if certain products of the exponents
$\rho_j$ are larger than one. Hence, our results are applicable not only to
the ac-networks, but to more general cases as well. The estimates
for local maps depend on the local structure of a considered network near
an equilibrium and are of a general form. We also state more cumbersome
conditions for asymptotic stability that take into account possible walks
that can be followed by a nearby trajectory.

Using these results we obtain conditions for asymptotic stability of two ac-networks and for two networks that do not belong to this class.
We also present  in Example~\ref{example5} two ac-networks whose stability was studied in \cite{afraimovich16}.
The networks in Example~\ref{example5} never satisfy our conditions for
asymptotic stability, while the conditions of \cite{afraimovich16} may
possibly be satisfied\footnote{For readers familiar with these examples or on a second reading of this article, we mention that in Example~\ref{example1}, if $|c_4/e_4|<1$ the conditions of \cite{afraimovich16}
are not satisfied, while ours may be.}.

The present article is organized as follows: the next section provides the necessary background for understanding our results. It is divided into three subsections, each focussing in turn on networks, graphs and stability. In Section~\ref{class}, we classify ac-networks. Section~\ref{sufcond} contains the main result providing sufficient conditions for the asymptotic stability of ac-networks. Section~\ref{exa} provides some illustrative examples.
The final section discusses the realisation of ac-networks and possible directions for continuation of the present study.

\section{Background and definitions}

\subsection{Robust heteroclinic networks}\label{sec:robust_network}

Consider a smooth dynamical system in $\R^n$ defined by
\begin{equation}
\dot{\bx} = f(\bx), \qquad \bx \in \R^n, \;\;
\label{eq:ode}
\end{equation}
and denote by $\Phi(\tau,\bx)$ the flow generated by solutions of the system.

We say the dynamical system is  $\Gamma$-equivariant, where $\Gamma$ is a compact Lie group, if
$$
 f(\gamma. \bx) = \gamma. f(\bx) \;\;\; \forall
\; \bx \in \R^n, \;\; \gamma \in \Gamma.
$$
Recall that the {\em isotropy subgroup} of $\bx$
is the subgoup of elements $\gamma\in\Gamma$ such that $\gamma\bx=\bx$. For a subgroup
$\Sigma\subset\Gamma$ the {\em fixed point subspace} of $\Sigma$ is
$$
\Fix(\Sigma)=\{~\bx\in\R^n~:~\sigma\bx=\bx~\hbox{ for all }\sigma\in\Sigma~\}.
$$
More detail about equivariant dynamical systems can be found in Golubitsky and Stewart \cite{SymmPerspective}.

For a hyperbolic equilibrium $\xi$ of \eqref{eq:ode} its stable and unstable
manifolds are
$$
W^u(\xi)=\{\bx\in\R^n~:~\alpha(\bx)=\xi\}, ~~~W^s(\xi)=\{\bx\in\R^n~:~\omega(\bx)=\xi\}.
$$
By $\kappa_{ij}$ we denote a
{\em connecting trajectory} from $\xi_i$ to $\xi_j$.
Following the terminology of Ashwin {\em et al} \cite{AshCasLoh18}, the
full set of connecting trajectories from $\xi_i$ to $\xi_j$,
$$
C_{ij}=W^u(\xi_i)\cap W^s(\xi_j),
$$
is called a {\em connection} from $\xi_i$ to $\xi_j$ and often also denoted by $[\xi_i\to\xi_j]$.
If $i=j$ then the connection is called {\em homoclinic}; otherwise, it is
{\em heteroclinic}.

A {\em heteroclinic cycle} is a union of a finite number of hyperbolic
equilibria, $\xi_1, \ldots ,\xi_m$, and heteroclinic connections $C_{j,j+1}$,
$1\le j\le m$, where $\xi_{m+1}=\xi_1$ is assumed.
A {\em heteroclinic network} is a connected union of finitely many
heteroclinic cycles.
In equivariant systems equilibria
and connections
 that belong to the same group orbit are usually identified.

A heteroclinic cycle is called \emph{robust} if every connection is of saddle-sink type in a flow-invariant subspace, $P_j$.
 In equivariant dynamical systems we typically have $P_j:=\Fix(\Sigma_j)$ for some subgroup $\Sigma_j\subset \Gamma$. In this case, due to $\Gamma$-invariance of fixed-point subspaces, the cycle persists with
respect to $\Gamma$-equivariant perturbations of $f$.

The eigenvalues of
the Jacobian, $\rd f(\xi_j)$,
with $\xi_j\in P_{j-1}\cap P_j$,
are divided into \emph{radial} (the associated
eigenvectors belong to $L_j=P_{j-1}\cap P_j$), \emph{contracting} (eigenvectors
belonging $P_{j-1}\ominus L_j$), \emph{expanding} (eigenvectors belonging
$P_j \ominus L_j$) and \emph{transverse} (the remaining ones), where
$P\ominus L$ denotes the complementary subspace of $L$ in $P$.

An equilibrium $\xi_j$ in a heteroclinic network might belong to several
heteroclinic cycles. In such a case the radial subspace is the same for all
cycles, since it is the subspace that is fixed by the isotropy subgroup of
$\xi_j$.
The stable manifold of $\xi_j$ contains all incoming connections at $\xi_j$.
The eigenvectors of $df(\xi_j)$ tangent to them are called \emph{contracting eigenvectors} at $\xi_j$
and they span the \emph{contracting subspace} at $\xi_j$.
The corresponding eigenvalues are called \emph{contracting eigenvalues.}
Similar notation is used for
the expanding and transverse eigenvalues, eigenvectors or subspaces.
By transverse eigenvalues we understand the eigenvalues that are not
radial, contracting or expanding for any cycle through $\xi_j$.

\begin{figure}[h]
\begin{center}
\parbox{0.45\textwidth}
{\begin{center}
\includegraphics[width=0.4\textwidth]{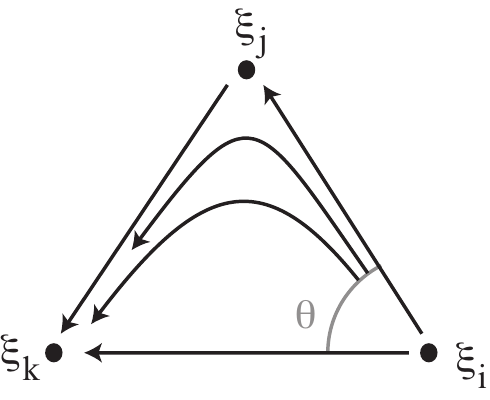}\\
(a)
\end{center}}
\qquad
\parbox{0.45\textwidth}
{\begin{center}
\includegraphics[width=0.4\textwidth]{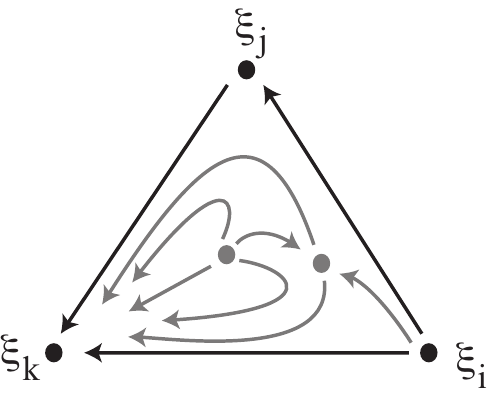}\\
(b)
\end{center}}
\end{center}
\caption{
(a) A $\Delta$-clique $\Delta_{ijk}$. (b)
Fragment of a phase portrait near a part of a network, where the union
of the connections $\xi_i\to\xi_j$, $\xi_i\to\xi_k$ and $\xi_j\to\xi_k$ is
not a $\Delta$-clique according to Definition \ref{Dcli}, but the respective
part of the associated graph is a $\Delta$-clique according to \cite{AshCasLoh18}.
 \label{figDeltaClique}}
\end{figure}

A heteroclinic network has subsets that are connected unions of
heteroclinic connections and equilibria and which are different from
heteroclinic cycles. An invariant set of the system (\ref{eq:ode})
which is called $\Delta$-clique is defined below. The name ``$\Delta$-clique''
was introduced in \cite{AshCasLoh18} to denote an element of a graph.
The relation between these two kinds of $\Delta$-cliques is discussed in the
following subsection, see also Figure~\ref{figDeltaClique}.

\begin{definition}\label{Dcli}
By a $\Delta$-clique, $\Delta_{ijk}$,
we denote the
union of equilibria and the following
 connecting trajectories  $\kappa_{\ell m}\subset C_{\ell m}=[\xi_\ell\to\xi_m]$:
\begin{equation}
\Delta_{ijk}
=\{\xi_i,\xi_j,\xi_k\}\cup\kappa_{ij}\cup\kappa_{jk}
\cup\biggl(\cup_{\theta\in[0,1)}\kappa_{ik}(\theta)\biggr),
\label{ddcli}
\end{equation}
such that
 $\kappa_{ik}(\theta)=\{\,\Phi(\tau,x_{ik}(\theta))\,:\,\tau\in\R\,\}$,
where $x_{ik}(\theta)$ is a continuous function of $\theta$,
$\kappa_{ik}(\theta_1)\ne\kappa_{ik}(\theta_2)$ for any $\theta_1\ne\theta_2$
and
$$\lim_{\theta\to1}\kappa_{ik}(\theta)=\kappa_{ij}\cup\xi_j\cup\kappa_{jk}.$$
\end{definition}

\subsection{Graphs and networks}

There is a close relation between heteroclinic networks and
directed graphs, which is discussed, e.g., by
Ashwin and Postlethwaite \cite{AshPos13}, Field \cite{fie15} or Podvigina
and Lohse \cite{pl19}. Given a heteroclinic network, there is an associated
graph such that the vertices (or nodes)
of the graph correspond to the equilibria of the network and an edge from $\xi_i$ to
$\xi_j$ corresponds to a full set of connections $C_{ij}$. (This supports the
frequent choice of the term {\em node} instead of equilibrium in the context of heteroclinic cycles and networks.)
Using the terminology from graph theory\footnote{We have used Foulds \cite{fou92} as a reference for graph theory.
There are many other good sources.}, we call a {\em walk} the union of vertices
$\xi_1,...,\xi_J$ and connections $C_{j,j+1}$, $1\le j\le J-1$.
The respective part of the graph is also called a {\em walk}. A walk where all
connections are distinct is called a {\em path}. A closed path
corresponds to a heteroclinic cycle, which is a subset of the network.

\begin{figure}[h]
 \centerline{\includegraphics[width=0.4\textwidth]{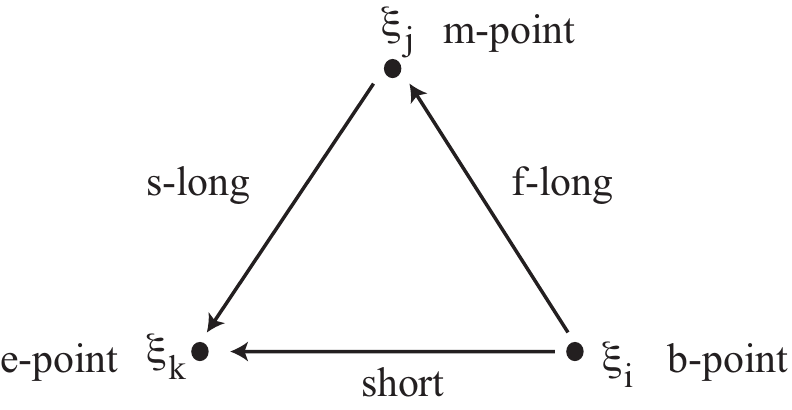}}
\caption{The $\Delta$-clique $\Delta_{ijk}$ with its short and long connections.\label{fig:Delta-clique}}
\end{figure}

In agreement with \cite{AshCasLoh18}, a nontransitive triangle within a graph
is called {\em $\Delta$-clique}: it is the union of three vertices $v_i$,
$v_j$ and $v_k$ and edges $[v_i \rightarrow v_k]$,
$[v_i \rightarrow v_j]$ and $[v_j \rightarrow v_k]$. In this case, in the
corresponding heteroclinic network the respective equilibria, $\xi_i$, $\xi_j$
and $\xi_k$, are connected by the trajectories
$\kappa_{ik}$, $\kappa_{ij}$ and $\kappa_{ik}$.
If a set $\Delta_{ijk}$ (see Definition \ref{Dcli}) is a subset
of a network, then the respective part of the graph is a $\Delta$-clique.
Hence, we use the term $\Delta$-clique to denote such a set.
Note that a network can possibly have a subset that is
represented by a $\Delta$-clique on the graph, but it is not
a $\Delta$-clique according to Definition \ref{Dcli},
since this subset does not necessarily contains
the full set of outgoing trajectories, as
required by Definition~\ref{Dcli}
(see Figure~\ref{figDeltaClique}~(b)).
However, this is not the case for the networks
considered in this paper.
We label equilibria and connections in a
$\Delta$-clique as follows
(the notation is shown in Figure~\ref{fig:Delta-clique}):

\begin{definition}\label{DDD}
Consider a $\Delta$-clique $\Delta_{ijk}$, which is an
invariant set of  system (\ref{eq:ode}).
We call the connections $[\xi_i\to\xi_j]$, $[\xi_j\to\xi_k]$ and $[\xi_i\to\xi_k]$,
the {\em f-long}, {\em s-long} and {\em short}  connections, respectively.
We call the equilibria $\xi_i$, $\xi_j$ and $\xi_k$, the {\em b-point},
{\em m-point} and {\em e-point}, respectively.
The eigenvectors of $\rd f(\xi_j)$ tangent to the
connecting trajectories
$\kappa_{ij}$ and $\kappa_{jk}$ are called {\em f-long} and
{\em s-long vectors}. Here there letters ``f'' and ``s'' stand for the first
and second; the letters ``b'', ``m'' and ``e'' stand for beginning, middle and
end, respectively.
\end{definition}

\subsection{Asymptotic stability}
\label{asstab}

Below we recall two definitions of asymptotic stability that can be found in
literature. In what follows $X\subset\R^n$ and $Y\subset\R^n$ are invariant
sets of the dynamical system.

\begin{definition}
A set $X$ is {\em Lyapunov stable} if for any neighbournood $U$ of $X$,
there exists a neighbourhood $V$ of $X$ such that
$$\Phi(\tau,\bx)\in U\hbox{ for all }\bx\in V\hbox{ and }\tau>0.$$
\end{definition}

\begin{definition}[Meiss \cite{Meiss}]\label{as}
A set $X$ is {\em asymptotically stable} if it is Lyapunov stable and
in addition the neighbourhood $V$ can be chosen such that
$$\lim_{\tau\to\infty}d(\Phi(\tau,\bx),X)\to0\hbox{ for all }\bx\in V.$$
\end{definition}

\begin{definition}[Oyama {\em et.al.} \cite{Oyama}]\label{asII}
A set $X$ is {\em asymptotically stable} if it is Lyapunov stable and
in addition the neighbournood $V$ can be chosen such that
$$\cup_{\bx\in V}\omega(\bx)\subseteq X.$$
\end{definition}

We refer to asymptotic stability according to Definition~\ref{as} as a.s. and
according to Definition~\ref{asII} as a.s.II. In this paper we adopt
Definition~\ref{as} as the definition of asymptotic stability.

\begin{figure}[h]
\begin{center}
\parbox{61mm}{
\begin{center}
\includegraphics[width=6cm]{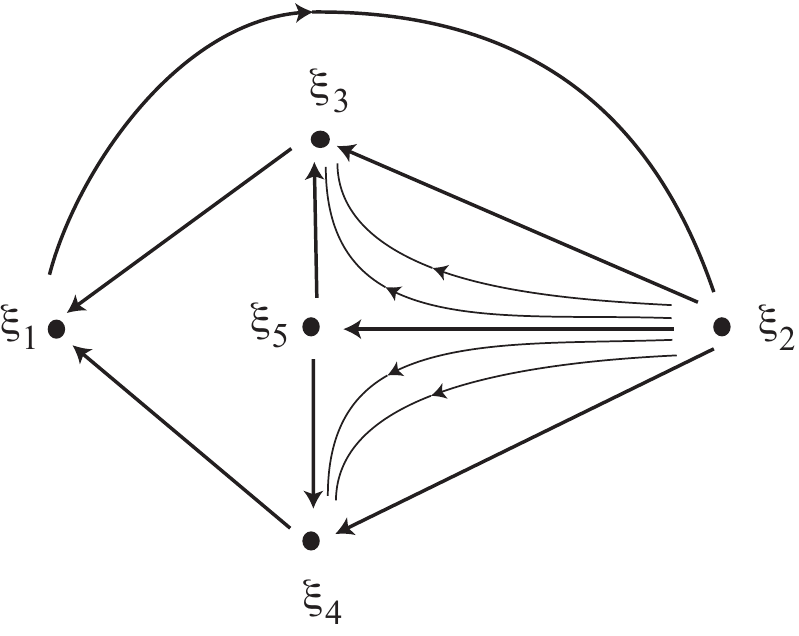}\\
(a)
\end{center}
}
\qquad
\parbox{61mm}{
\begin{center}
\includegraphics[width=6cm]{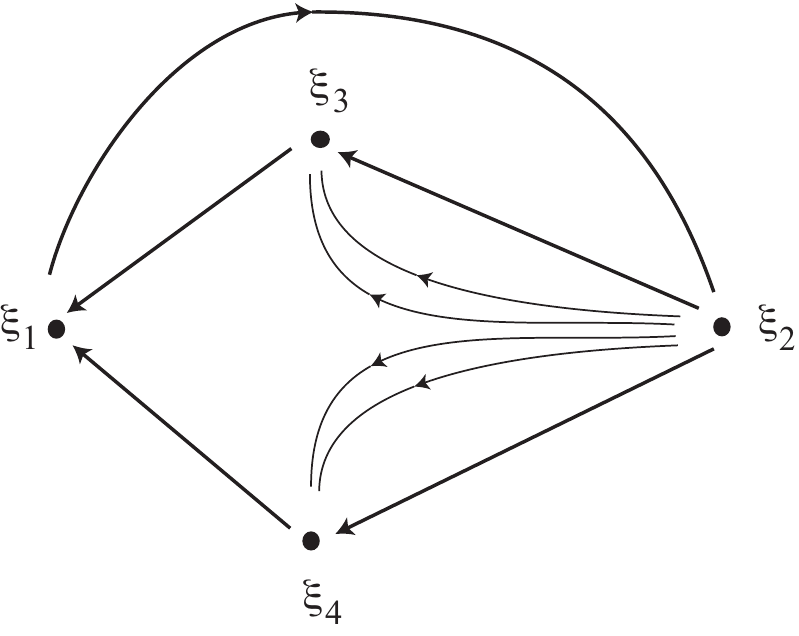}\\
(b)
\end{center}
}
\end{center}

\caption{
A clean heteroclinic network (a) and  (b) its subset $X$, which is not clean.
The network in (a) can possibly be a.s.II, while the one
in (b) is at most a.s. because there exist points $x$ in a neighbourhood of $\xi_2$
such that  $\omega(x)=\xi_5\notin X$.
}\label{fig2}
\end{figure}

\begin{remark}\label{rem30}
For closed sets, such as equilibria or periodic orbits, the Definitions~\ref{as}
and \ref{asII} are equivalent. However, this is not always the case.
In Figure \ref{fig2}b we show an example of
a heteroclinic network, which is not a.s.II, but is a.s.
if the eigenvalues at equilibria satisfy the conditions stated
in Section~\ref{exa} (see Example~\ref{example3}).
Note also the following relations between these two definitions of asymptotic stability:
$$\hbox{($X$ is a.s.II)$\Rightarrow$($X$ is a.s.)},~~~~
\hbox{($X$ is a.s.)$\Rightarrow$($\overline{X}$ is a.s.II)},~~~~
\hbox{($\overline{X}$ is a.s.II)$\Rightarrow$($X$ is a.s.)}$$
\end{remark}

The next definition is a generalization of Definition 1.3 of Field~\cite{fie15}, where  a  heteroclinic network is defined as {\em clean} if it is compact and equal to the union of the unstable manifolds
of its nodes.

\begin{definition}\label{clean}
(Adapted from \cite{fie15}.) A compact invariant set $Y$ is {\em clean}
if each of its invariant subsets $X$ satisfies
$$W^u(X)\subseteq Y.$$
\end{definition}

In what follows we denote by $N_\varepsilon(Y)$ the set $N_\varepsilon(Y)=\{\bx: d(\bx,Y)<\varepsilon \} $ for a closed set $Y$.

\begin{theorem}\label{th0}
If a compact invariant set is asymptotically stable then it is clean.
\end{theorem}

\proof
We show that if a compact invariant set $Y$ is not clean then it is not Lyapunov stable.

If $Y$ is not clean, there exists an invariant set $X \subset Y$ and $\bx\in W^u(X)$
such that $\bx \notin Y$. Since $Y$ is a compact set, we know that $h=d(\bx,Y)>0$.
Let $U= N_{h/2}(Y)$ and for any $V \supset Y$ take $\delta > 0$ such that
$N_\delta (X) \subset V$. Since $\bx\in W^u(X)$,
there exists at least one $\bx_\delta \in \Phi(\tau,\bx) \cap N_\delta(X) \subset V$.
It is clear that solutions through $\bx$ and $\bx_\delta$ coincide so that there
exists $\tau_\delta$ such that $\bx=\Phi(\tau_\delta, \bx_\delta) \notin U$.
\qed

\begin{Cy}\label{cleannetw}
Suppose that the system (\ref{eq:ode}) has a compact heteroclinic network,
which is a union of nodes (they can be any invariant sets) and heteroclinic
connections. If the network is a.s. then it is clean.
\end{Cy}

\section{Classification of ac-networks}
\label{class}

In this section we define ac-networks  and classify them by describing all possible types of associated graphs.
The letters ``ac'' stand for {\em axial} and {\em clean}.

\begin{definition}\label{acn}
An {\em ac-network} is a robust heteroclinic network $Y \subset \R^n_+$ in a $\Z_2^n$-equivariant system \eqref{eq:ode},  such that all the equilibria are hyperbolic  and lie on coordinate axes, there is no more that one equilibrium per axis, the dimension of the unstable manifold of any equilibrium is either one or two, and the network is clean.
\end{definition}

For every equilibrium $\xi_j$ in an ac-network, the $\Z_2^n$-equivariance of the vector field ensures that the eigenvalues of $\rd f(\xi_j)$ are real and the eigenvectors are aligned with the Cartesian basis vectors.
By application of the $\Z_2^n$ symmetries, the network may be extended from $\R^n_+$ to  $\R^n$.

On the other hand, the $\Z_2^n$ symmetries prevent some well known networks (whose graph may resemble that of an ac-network) from being ac-networks.
As an illustration, we refer to the network in Postlethwaite and Dawes \cite{PD2005} which has symmetry $\Z_3\ltimes Z_2^6$ on $\R^3\oplus\R^3$, or the Guckenheimer and Holmes cycle \cite{gh88} in $\R^3$, which also has symmetries other than $\Z_2^n$.

The demand that all equilibria are on coordinate axes excludes the Ashwin \emph{et al} \cite{AshCasLoh18} completion of the Kirk and Silber network from the set of ac-networks.
Note that an ac-network may exist for a vector field possessing equilibria outside the coordinate axes.
The demand in Definition~\ref{acn} concerns only the equilibria that are part of the network.

The structure of ac-networks is natural in population dynamics, see examples in Afraimovich \emph{et al} \cite{afraimovich16}.
The simplex method of Ashwin and Postlethwaite \cite{AshPos13} produces networks that are axial and have $\Z_2^n$ symmetry.
They are not necessarily clean.

\bigbreak
In the following two lemmas $Y$ denotes an ac-network.
Recall that $\kappa_{ij}$ is a connecting trajectory in full set of connections from $\xi_i$ to $\xi_j$ denoted by $C_{ij}=[\xi_i\to \xi_j]$.

\begin{lemma}\label{lem0}
An ac-network
does not have subsets that are homoclinic cycles or
heteroclinic cycles with two equilibria.
\end{lemma}

\proof
By Definition \ref{acn} any connection $C_{ij} =W^u(\xi_i)\cap W^s(\xi_j) \subset Y$ is
robust, i.e., it belongs to a flow-invariant subspace where $\xi_i$ is
unstable and $\xi_j$ is a sink. Since an unstable equilibrium is not a
sink, an ac-network does not have homoclinic connections.

Suppose there exists a heteroclinic cycle with two equilibria $\xi_i$ and $\xi_j$ and connections $C_{ij}$ and $C_{ji}$. Let $S_{ij}$ be the flow-invariant subspace containing $C_{ij}$. In the coordinate plane $P_{ij}=<\be_i,\be_j> \subset S_{ij}$ the trajectory $\kappa_{ij} = P_{ij}\cap C_{ij}$ connects $\xi_i$ to $\xi_j$. Similarly, there exists a trajectory $\kappa_{ji} = P_{ij} \cap C_{ji}$ connecting $\xi_j$ to $\xi_i$. Since $\xi_i$ and $\xi_j$ are hyperbolic equilibria, the connecting trajectories $\kappa_{ij}$ and $\kappa_{ji}$ cannot exist robustly and simultaneously in $P_{ij}$.
\qed

\begin{lemma}\label{lem1}
An ac-network
is a union of equilibria, one-dimensional connections
and $\Delta$-cliques.
\end{lemma}

\proof
Definition~\ref{acn} implies that all connections are 1- or 2-dimensional.
Hence, due to the compactness of the network and the definition of $\Delta$-clique it suffices to show that the closure of any 2-dimensional connection is a $\Delta$-clique.

\begin{figure}[h]
 \centerline{\includegraphics[scale=.85]{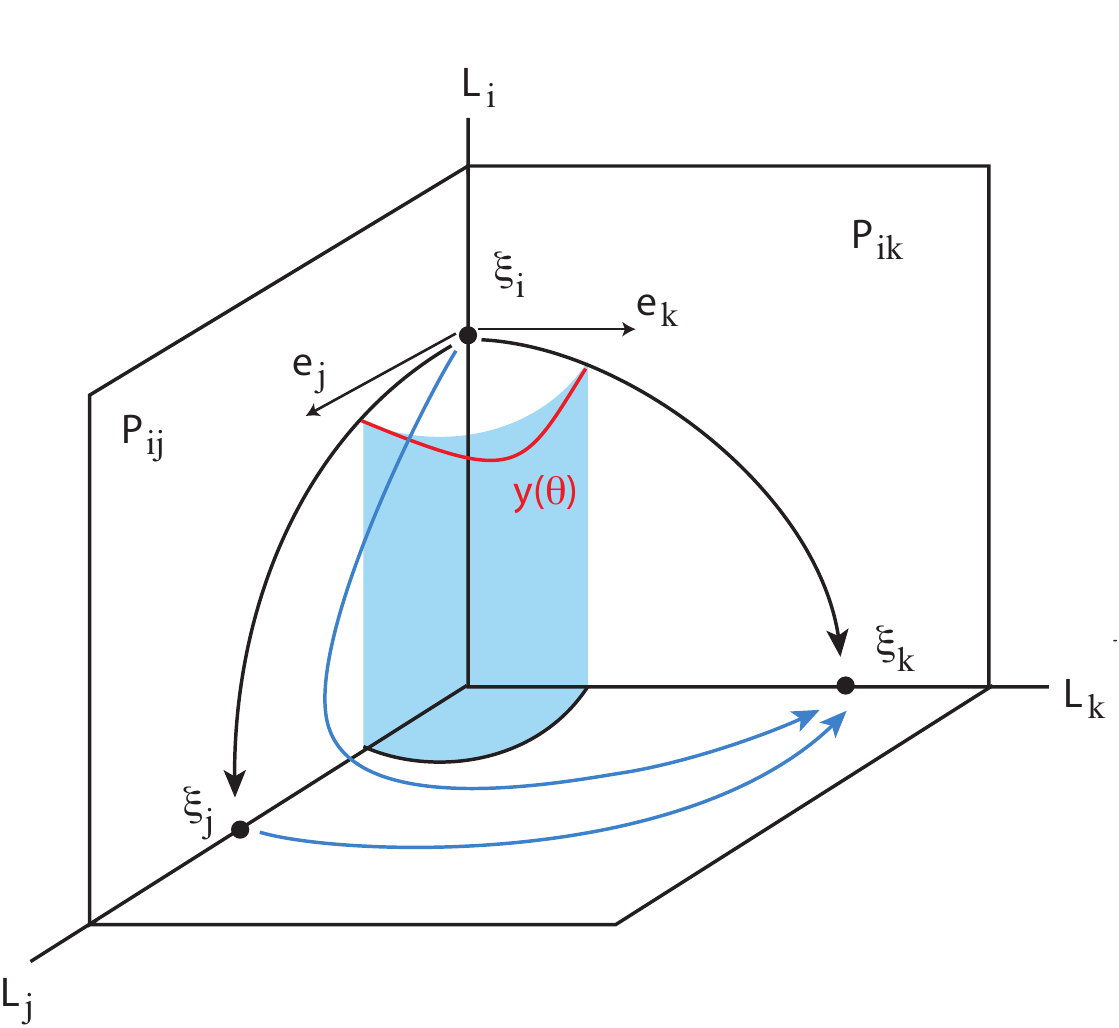}}
\caption{Conventions in the proof of Lemma~\ref{lem1} (case (i))
where $df(\xi_i)$ has expanding eigenvectors $\be_j$ and $\be_k$.
Trajectories through the curve $y(\theta)$ where $W^u(\xi_i)$ meets a small cylinder (shaded) around the $L_i$ axis must approach either $\xi_j$ or $\xi_k$. Then all trajectories through the part of $y(\theta)$ that lie outside the planes $P_{ij}$ and $P_{ik}$ must go to the same equilibrium (here shown to be $\xi_k$). Otherwise there would be a point $y(\theta_1)$ separating $W^s(\xi_j)\cap \{y(\theta)\}$ from $W^s(\xi_k)\cap \{y(\theta)\}$ such that $\omega(y(\theta_1))\notin \{\xi_j,\xi_k\}$.
\label{fig:Lemma33}}
\end{figure}

Let $\dim(W^u(\xi_i))=2$ and the expanding eigenvectors be $\be_j\subset L_j$
and $\be_k\subset L_k$. The planes  $P_{ij}=<\be_i,\be_j>$ and
$P_{ik}=<\be_i,\be_k>$ are flow-invariant, which implies the existence of
 the steady states
$\xi_j=\omega(W^u(\xi_i)\cap P_{ij})$ and
$\xi_k=\omega(W^u(\xi_i)\cap P_{ik})$
and the connecting trajectories $\kappa_{ij} \subset P_{ij}$ and $\kappa_{ik} \subset P_{ik}$
(see Figure~\ref{fig:Lemma33}).

Consider the dynamics within the invariant subspace
$S_{ijk}=<\be_i,\be_j,\be_k>$. Define a curve $\by(\theta)$,
$\theta\in[0,1]$ as the intersection of $W^u(\xi_i)$ with the surface described by
$(x_j,x_k)=(\delta\sin(\pi\theta/2),\delta\cos(\pi\theta/2))$,
where $\delta$ is a small number.
The curve can be written as
$\by(\theta)=(x_i(\theta),\delta\sin(\pi\theta/2),\delta(\cos\pi\theta/2))$.
A trajectory through $\by(\theta)$ satisfies
$$
\alpha(\Phi(\tau,\by(\theta)))=\xi_i\hbox{ and }
\omega(\Phi(\tau,\by(\theta)))=\xi_j\hbox{ or }\xi_k.
$$
Because the equilibria of an ac-network belong to coordinate axes
 the network contains no other equilibrium in $S_{ijk}$.
Since the network is clean,
there can be no $\theta_1$ such that $\omega(\Phi(\tau,\by(\theta_1)))=\xi^\prime \ne \xi_j,\xi_k$
 as $\xi^\prime$ would be in $S_{ijk}$.
Suppose for some $\theta_1\in(0,1)$ and some $\varepsilon>0$ we had $\omega(y(\theta))=\xi_j$ for $ \theta_1-\varepsilon<\theta<\theta_1$ and
$\omega(y(\theta))=\xi_k$ for $ \theta_1<\theta<\theta_1+\varepsilon$.
Then $\omega(y(\theta_1))$ could not be contained in the network, contradicting the fact that it is clean.
Then one of the following two alternatives holds:
\begin{itemize}
\item [(i)] $\omega(\Phi(\tau,\by(1)))=\xi_j$ and for $0\le \theta<1$ we have $\omega(\Phi(\tau,\by(\theta)))=\xi_k$;
\item [(ii)] $\omega(\Phi(\tau,\by(0)))=\xi_k$ and for $0<\theta\le 1$ we have $\omega(\Phi(\tau,\by(\theta)))=\xi_j$.
\end{itemize}
In case (i) due to the compactness of $Y$ and the smoothness of $f$
$$
\lim_{\theta\to1}\Phi(\tau,\by(\theta))=\left(\xi_j\cup\kappa_{ij}\cup\kappa_{jk}\right)\subset Y.
$$
Therefore, the $\Delta$-clique $\Delta_{ijk}$
is a subset of $Y$. In case (ii) there exists
a $\Delta$-clique $\Delta_{ikj}$ which is a subset of $Y$.
\qed

\begin{Cy}\label{outcon}
If $C_{ij},C_{ik}\in Y$ then either $C_{jk}\subset Y$ or $C_{kj}\subset Y$. Furthermore, if $C_{ik}$ is not
an f-long connection then $C_{jk}\subset Y$ hence $\Delta_{ijk}\subset Y$.
\end{Cy}

We claim that an ac-network always contains at least one cycle without f-long connections.
Indeed, if one starts from any node, either it has only one  outgoing connection (hence not f-long), or it has one short and one f-long. Follow the short (or only) connection to the next node and repeat the reasoning there. Since the set of nodes is finite, this eventually closes into a cycle. See also Remark~\ref{endcycl}.
The claim implies that the following theorem describes all possible ac-networks.

\begin{theorem}\label{th1}
Let $Y\in\R^n$ be an ac-network and let
$X=(\cup_{1\le j\le J}\xi_j)\cup(\cup_{1\le j\le J}C_{j,j+1})$ be a
heteroclinic cycle in $Y$. Assume that $X$ has no f-long connections.
Then either
\begin{itemize}
\item [(i)]  there is at least one connection in $Y$ between equilibria in $X$ which is not contained in $X$.
Then $J$ is odd and $Y$ is a union of equilibria $\xi_j$,
$1\le j\le J$, connected by $\Delta$-cliques.
\item [(ii)]  all connections in $Y$ between equilibria in $X$ are contained in $X$.
Then the equilibria $\xi\in Y$ can be grouped into disjoint sets
$$
\cX_1=\{\xi_{1,0},...,\xi_{1,m_1}\},...,\cX_J=\{\xi_{J,0},...,\xi_{J,m_J}\}
$$
where $\xi_{j,0}=\xi_j\in X$.
The network $Y$ always has the connections:
$$
[\xi_{j,s}\to\xi_{j+1,0}],\quad 0\le s\le m_j,\quad
[\xi_{j,s-1}\to\xi_{j,s}],\quad 1\le s \le m_j,
$$
where $j=1, \hdots, J$.
The network $Y$ may also have one of the following connections
$$
[\xi_{j,m_j}\to\xi_{j+1,1}], \qquad
[\xi_{j,m_j}\to\xi_{j+2,0}],\qquad [\xi_{j,m_j}\to\xi_{j,0}].
$$
It may also have  the connection $[\xi_{j-1,m_{j-1}}\to\xi_{j,m_j}]$,  if the connection $[\xi_{j,m_j}\to\xi_{j,0}]$ exists.
\end{itemize}
\end{theorem}

\proof
(i) Let $C_{1s}=\,[\xi_1\to\xi_s]\subset Y$, $C_{1s}\not\subset X$, be
the connection in the statement of the theorem. The equilibrium $\xi_1$
has two outgoing connections, $C_{12}\subset X$ and $C_{1s}$. Therefore
by Corollary \ref{outcon}, since $C_{12}$ is not f-long, there exists a further connection $C_{s2}\subset Y$ and the corresponding $\Delta$-clique, $\Delta_{1s2}$.

Hence, the equilibrium $\xi_s$ has two outgoing connections, $C_{s,s+1}$
and $C_{s2}$. By the same arguments as above, this implies existence
of the connection $C_{2,s+1}$, which in turn implies existence of the
connection $C_{s+1,3}$. Proceeding like this, we prove existence of the
connections $C_{3,s+2}$, $C_{4,s+3}$, and so on up to
$[\xi_s\to\xi_{2s-1}]$. However,
$\xi_s$ already has two outgoing connections, $[\xi_s\to\xi_{s+1}]$
and $[\xi_s\to\xi_2]$. An equilibrium in an ac-network has one or two
outgoing connections, therefore $2s-1 \pmod{J}=2$, which implies that
$2s-3=J$.

Note, that we have identified other
$J$ connections $[\xi_i\to\xi_{i+(J+1)/2}]$,
each of them is an f-long connection for one $\Delta$-clique and s-long
connection for another. A connection $C_{j,j+1}\subset X$
is the short connection of two $\Delta$-cliques.
Figure~\ref{fig:Ex5} (left panel) shows examples of this type of network.

\bigskip
(ii) Next we consider the case when $Y$ does not have connections between
$\xi_j$, $1\le j\le J$, other than $[\xi_j\to\xi_{j+1}]$.
To each equilibrium $\xi_j$ in $X$ we associate a class $\cX_j$.

The grouping of the equilibria in $Y$ into the sets $\cX_j$ proceeds stepwise: after assigning each equilibrium $\xi_j$ in $X$ to a group $\cX_j$ (step 1), subsequent steps distribute the equilibria in $Y\backslash X$ by each $\cX_j$, or move them from one group to another, in such a way that the statement of the theorem is satisfied.
In the construction we use
an auxiliary network $Z\subseteq Y$, which is modified on each step by adding
equilibria (in parallel with $\cX$)  and/or connections. At Step 1, $Z=X$ and, because the number of equilibria in $Y$ is finite, this process is finite and at the last step $Z=Y$.

It follows from the assumptions of case (ii) of the theorem that there is at least one equilibrium $\xi \in Y$ but $\xi \notin X$, and a connection from some $\xi_j \in X$ to $\xi$ so that Step 2 is always taken. There are at least two connections of the compulsory type listed in the statement of the theorem.
Step 3 accounts for possible connections between the equilibria already classified in Step 2.
Steps 4 and 5 are required only if there are $\xi_{.,1}$ with connections to other equilibria not already in $Z$.

 The number of connections between equilibria is restricted by the hypotheses and by Corollary~\ref{outcon}.
\bigbreak

\paragraph{Step 1: initiation.} We initiate the sets $\cX_j=\{\xi_{j,0}\}$, $1\le j\le J$, where
$[\xi_{j,0}=\xi_j]$ are the equilibria of $X$ and define $Z=X$.

\paragraph{Step 2: compulsory connections for $s \le 1$.} Some of $\xi_{j,0}=\xi_j$ have outgoing connections other
than $[\xi_j\to \xi_{j+1}]$, i.e., connections $C_{jk}\subset Y$ where
$k>J$. For each $j$ there is at most one such connection due to the definition
of ac-networks. Existence of a connection $C_{jk}$ implies
that there are no connections $C_{ik}$ from other $\xi_i\in X$,
$i\ne j$, to $\xi_k$.\footnote
{Since $C_{j,j+1}$ is not an f-long connection for any $\Delta$-clique,
by Corollary \ref{outcon}, existence of the connections $C_{jk}$
and $C_{j,j+1}$ implies the existence of the connection
$C_{k,j+1}$. Similarly, existence of the connections
$C_{ik}$ and $C_{i,i+1}$ implies the
existence of $C_{k,i+1}$. Hence, the corollary implies existence of a
connection between $C_{j+1}$ and $C_{i+1}$, which contradicts the statement
of the theorem unless $i=j+1$ or $j=i+1$. But in the case, e.g. $i=j+1$, there exist
both connections $C_{j+1,k}$ and $C_{k,j+1}$, which is not
possible in an ac-network.}
So, if a connection $[\xi_j\to\xi_k]$ exists,
we assign $\xi_k$ to $\cX_j$: $\xi_{j,1}=\xi_k$. We include in $Z$ all such equilibria,
the connections $[\xi_j\to\xi_k]$ and $[\xi_k\to\xi_{j+1}]$. The latter exists
due to Corollary \ref{outcon} since the connection $[\xi_j\to\xi_{j+1}]$
is not an f-long connection of a $\Delta$-clique. Note, that at the end of
this step all outgoing connections from $\xi_{j,0}$ belong to $Z$.

\begin{figure}[h]

\vspace*{5mm}
\centerline{
\includegraphics[scale=.70]{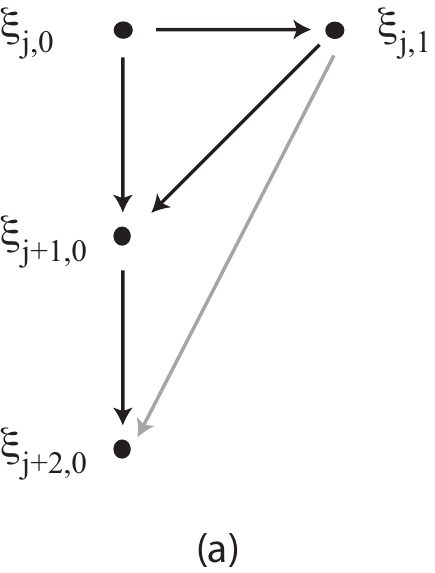}\quad
\includegraphics[scale=.70]{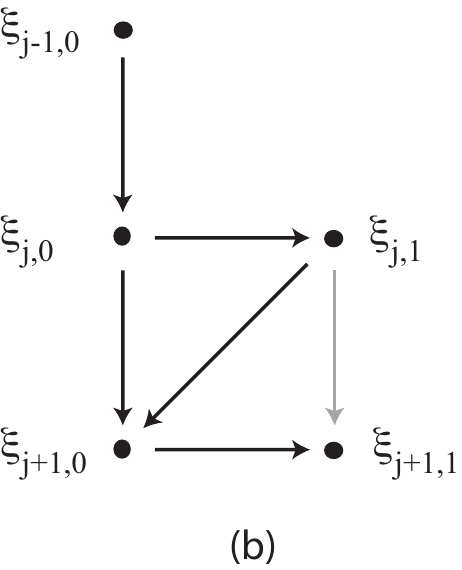}\quad
\includegraphics[scale=.70]{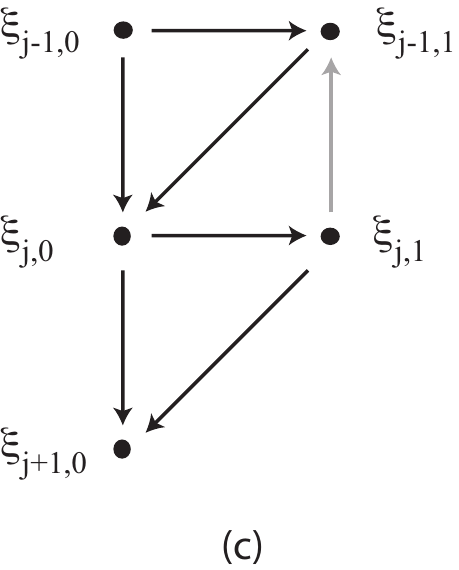}\quad
\includegraphics[scale=.70]{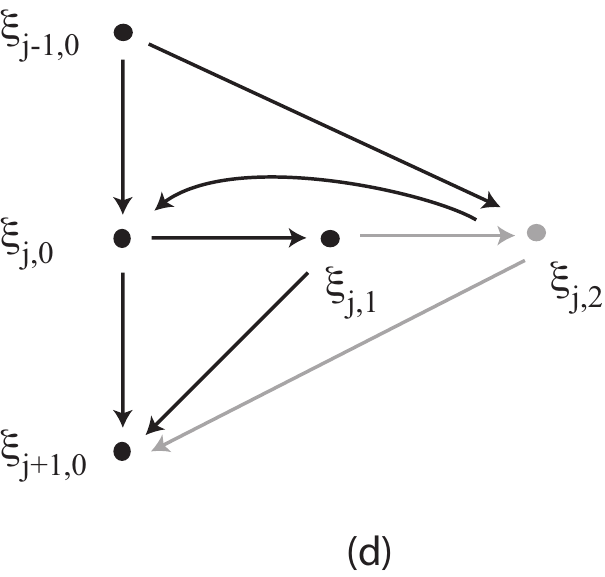}
}

\caption{Possible cases for optional connections in Step 3 of  the construction of the network $Z$, new optional connections and consequent compulsory connections shown in grey. The structure in (c) is  moved to that in (d) generating another compulsory connection. In all these examples the set $\cX_j$ is complete at this step.}\label{figStep3}
\end{figure}

\paragraph{Step 3: optional connections for $m_j=1$.} We add to $Z$ connections $[\xi_{j,1}\to \xi_i]$, such that
$\xi_i\in Z$ and the connection was not included in $Z$ at the previous step.
These are of two types: either from $\xi_{j,1}$ back to an equilibrium in $X$ or from $\xi_{j,1}$ to an equilibrium $\xi_{s,1}$.

By Lemma \ref{lem1}, if a connection from $\xi_{j,1}$ back to an equilibrium in $X$ exists, it must be
$[\xi_{j,1}\to\xi_{j+2,0}]$.
In this case $m_j=1$, since $\xi_{j,1}$ already has two outgoing connections,
$[\xi_{j,1}\to\xi_{j+1,0}]$ (identified in Step 2) and $[\xi_{j,1}\to\xi_{j+2,0}]$
as in Figure~\ref{figStep3} (a).

Lemma~\ref{lem1} gives two possibilities for
$[\xi_{j,1}\to\xi_{s,1}]$: they are $[\xi_{j,1}\to\xi_{j+1,1}]$ or
$[\xi_{j,1}\to\xi_{j-1,1}]$. (There are no other possibilities due to
Corollary~\ref{outcon} and absence in $Y$ of
additional connections between $\xi_j\in X$.)
In the former case
(Figure~\ref{figStep3} (b))
 we have $m_j=1$.
In the latter case
(Figure~\ref{figStep3} (c))
we move $\xi_{j-1,1}$ from $\cX_{j-1}$ to $\cX_j$ where we set
$\xi_{j,2}=\xi_{j-1,1}$
(Figure~\ref{figStep3} (d)).
Then there are two outgoing connections from $\xi_{j,1}$ and by Corollary~\ref{outcon}
there is a connection between $\xi_{j,2}$ and $\xi_{j+1,0}$.
The connection   $[\xi_{j+1,0}\to \xi_{j,2}]$  is not compatible with the existing $[\xi_{j-1,0}\to \xi_{j,2}]$, as shown in Step 2.
Hence, the second connection is $[\xi_{j,2}\to \xi_{j+1,0}]$, that we add to  $Z$, and $m_j=2$, $m_{j-1}=0$.

\begin{figure}[h]

\vspace*{5mm}
\centerline{
\includegraphics[scale=.70]{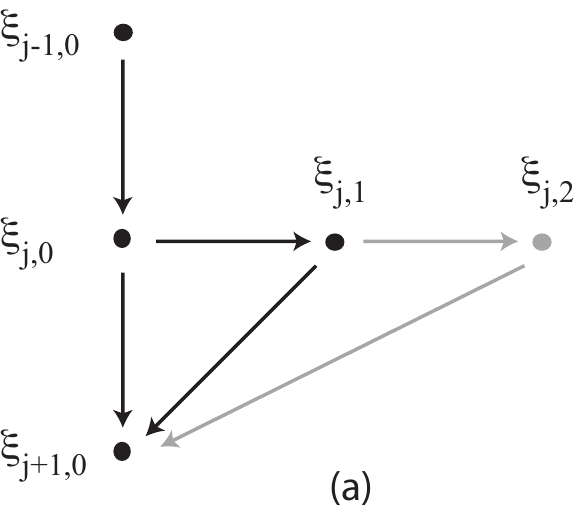}\qquad
\includegraphics[scale=.70]{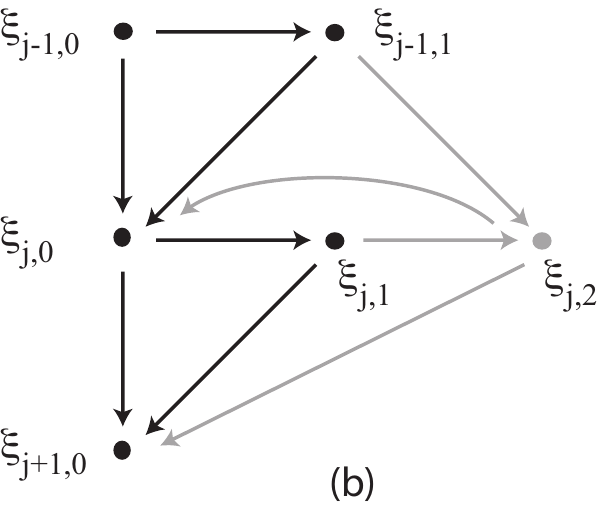}
}

\caption{Compulsory (a) and  optional (b) connections in Step 4 of  the construction of the network $Z$, new connections and consequent compulsory connections shown in grey. In case (b) the set $\cX_j$ is complete at this step.}\label{figStep4}
\end{figure}

\paragraph{Step 4: compulsory connections for $s=2$.} Let $\xi_{j,1}$ have an outgoing connection $[\xi_{j,1}\to\xi_i]$
where $\xi_i\not\in Z$ and assume there are no other connections $[\xi_{k,1}\to\xi_i]$.
Then we add $\xi_i$ to $\cX_j$: $\xi_{j,2}=\xi_i$.
As well, we add to $Z$ the equilibrium
$\xi_{j,2}=\xi_i$ and the connections
$[\xi_{j,1}\to\xi_i]$ and $[\xi_i\to\xi_{j+1}]$
(Figure~\ref{figStep4} (a)).
If some $\xi_i\in Y$ has several such incoming connections from the
equilibria $\xi_{j,1}$, then there are two such connections and they are
$[\xi_{j-1,1}\to\xi_i]$ and $[\xi_{j,1}\to\xi_i]$. (This can be shown by arguments
similar to the ones applied at Step 2.)
Then we add $\xi_i$ to $\cX_j$: $\xi_{j,2}=\xi_i$
(Figure~\ref{figStep4} (b)).
 We add to $Z$
the equilibrium $\xi_i$ and the connections $[\xi_{j-1,1}\to\xi_i]$,
$[\xi_i\to\xi_j]$, $[\xi_{j,1}\to\xi_i]$ and $[\xi_i\to\xi_{j+1}]$.
Note that at the end of this step all outgoing
connections from $\xi_{j,1}$ that are in $Y$ belong to $Z$.

\begin{figure}[h]

\vspace*{5mm}
\centerline{
\includegraphics[scale=.70]{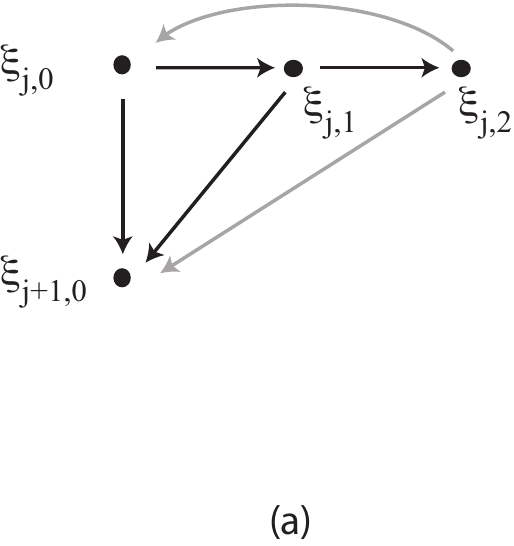}\quad
\includegraphics[scale=.70]{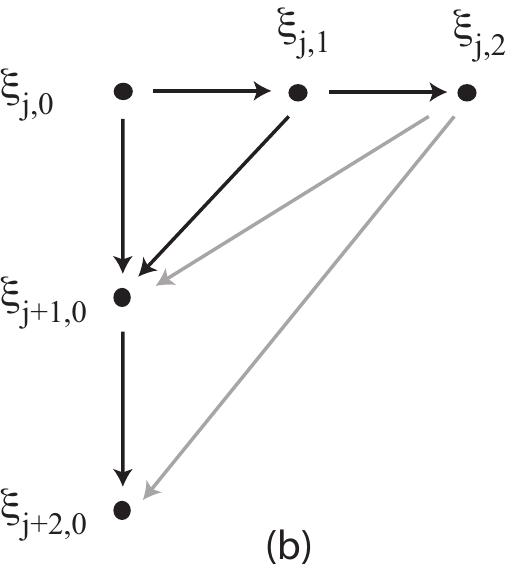}\quad
\includegraphics[scale=.70]{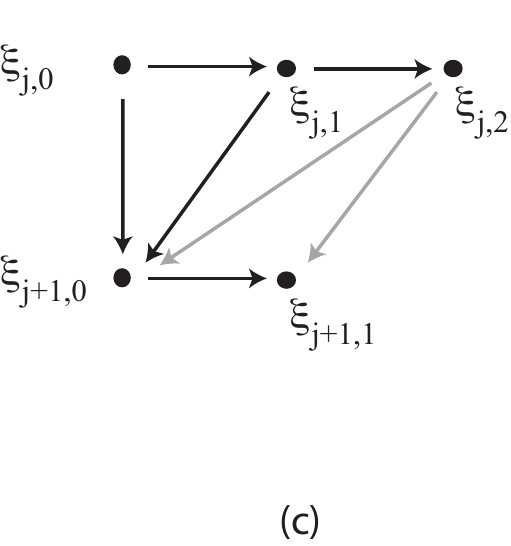}
}
\caption{Some possible cases for optional connections in Step 5 of  the construction of the network $Z$, new connections and consequent compulsory connections shown in grey. In all these examples the set $\cX_j$ is complete at this step. }\label{figStep5}
\end{figure}

\begin{figure}[h]

\vspace*{5mm}
\centerline{
\includegraphics[scale=.65]{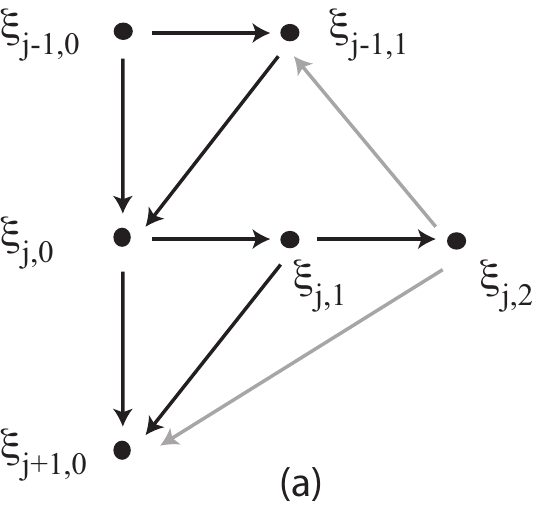}\quad
\includegraphics[scale=.65]{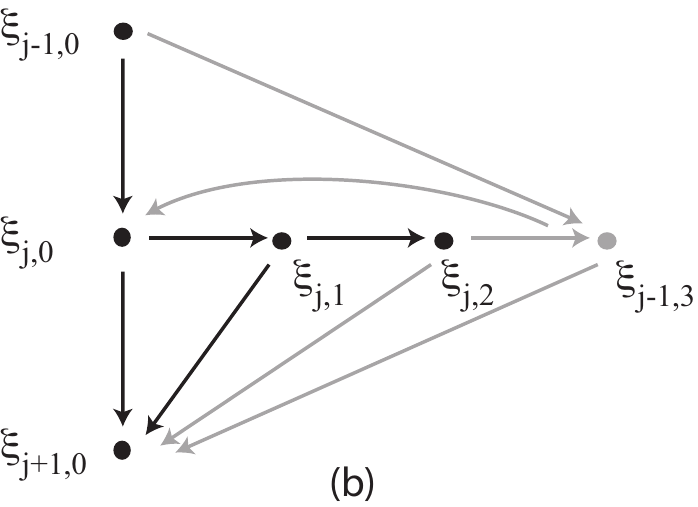}\quad
\includegraphics[scale=.65]{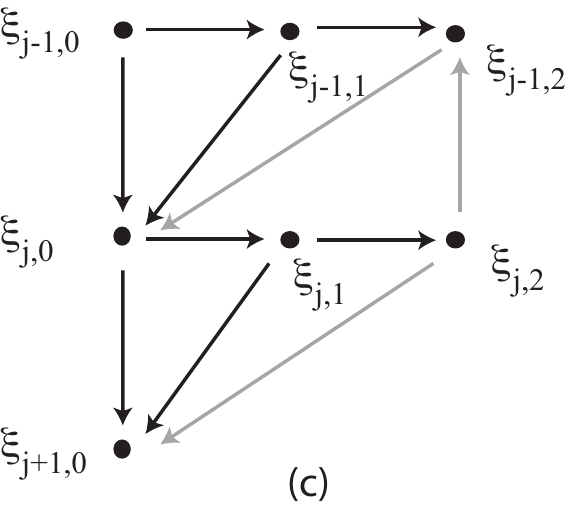}\quad
\includegraphics[scale=.65]{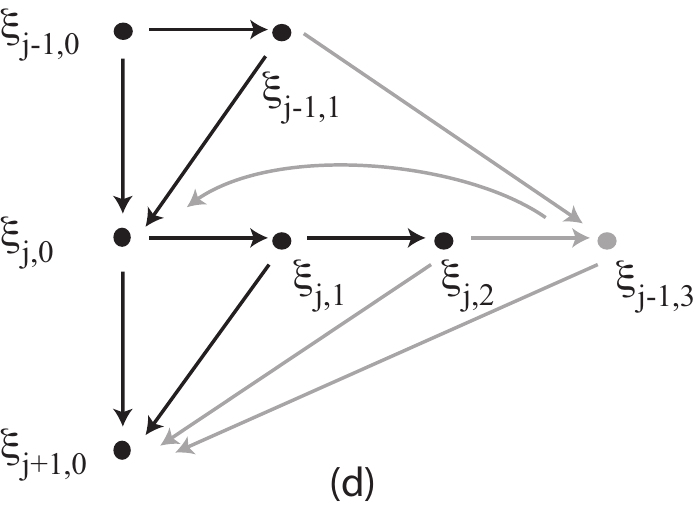}
}
\caption{Some possible cases for optional connections in Step 5 of  the construction of the network $Z$, where some $\xi_i$ is moved from $\cX_{j-1}$ to $\cX_j$.
New connections and consequent compulsory connections shown in grey.The structure in (a) is  moved to that in (b) and contains a compulsory connection. Similarly, (c) is moved to (d) In all these examples the set $\cX_j$ is complete at this step. }\label{figStep5A}
\end{figure}

\paragraph{Step 5: optional connections for $m_j=2$.} We add to $Z$ outgoing
connections from $\xi_{j,2}$ to
other equilibria that are already in $Z$, namely
$\xi_{s,0}$, $\xi_{s,1}$ and $\xi_{s,2}$.
Lemma \ref{lem1} implies that the only (mutually exclusive) possibilities are:
$[\xi_{j,2}\to\xi_{j,0}]$, $[\xi_{j,2}\to\xi_{j+2,0}]$,
$[\xi_{j,2}\to\xi_{j+1,1}]$, $[\xi_{j,2}\to\xi_{j-1,1}]$ or
$[\xi_{j,2}\to\xi_{j-1,2}]$
(Figure~\ref{figStep5} (a)--(c) and Figure~\ref{figStep5A} (a),(c)).
If a connection $[\xi_{j,2}\to\xi_{j-1,1}]$ exists
(Figure~\ref{figStep5A} (a)),
then
we move $\xi_{j-1,1}$ from $\cX_{j-1}$ to $\cX_j$ where we set
$\xi_{j,3}=\xi_{j-1,1}$
(Figure~\ref{figStep5A} (b)).
 If a connection $[\xi_{j,2}\to\xi_{j-1,2}]$ exists
(Figure~\ref{figStep5A} (c)),
then we move $\xi_{j-1,2}$ from $\cX_{j-1}$ to
$\cX_j$ where we set $\xi_{j,3}=\xi_{j-1,2}$
(Figure~\ref{figStep5A} (d)).
If one of the first three connections exists
then $m_j=2$, otherwise $m_j=3$.
\medbreak

Steps 4 and 5 are repeated
with larger $s$ until
 all equilibria of $Y$ become equilibria
of $Z$. Since $Y$ has a finite number of equilibria the process will be
terminated after a finite number of steps.

Therefore, we have obtained the network $Z$ that coincides with $Y$.
By construction, the network $Z$ has the structure stated in the theorem.
\qed
\bigbreak

The  decomposition of the sets of nodes in Theorem~\ref{th1} implies that an ac-network $Y$ contains at most two heteroclinic cycles without f-long connections.

\begin{remark}\label{endcycl}
The structure of the graphs associated with ac-networks implies that for a given ac-network  $Y$, the cycle $X$ without f-long connections of Theorem~\ref{th1} is  unique, except for the case when $Y$ is  comprised of equilibria
$\xi_{j,0}$ and $\xi_{j,1}$ and the connections
$$[\xi_{j,0}\to\xi_{{j+1},0}],\quad [\xi_{j,0}\to\xi_{j,1}],\quad
[\xi_{j,1}\to\xi_{{j+1},0}],\quad [\xi_{j,1}\to\xi_{j+1,1}],$$
where $1\le j\le J$.
For such $Y$ the cycle $[\xi_{1,1}\to\xi_{2,1}\to...\to\xi_{J,1}\to\xi_{1,1}]$
also has no f-long connections. See the right panel of Figure~\ref{fig:Ex5}.
\end{remark}

The sets $\cX_j$ in Theorem~\ref{th1} depend on the cycle $X$, so when  the network contains two cycles  without f-long connections this yields two different groupings.

\begin{remark}\label{cycleKS}
Recall the network in Figure~\ref{figKS} (b). Relabelling the equilibria as
$\zeta_{1,0}=\xi_2$;  $\zeta_{2,0}=\xi_3$;  $\zeta_{3,0}=\xi_1$;
 $\zeta_{1,1}=\xi_4$, the compulsory connections are
$[\zeta_{s,0}\to \zeta_{s+1,0}]$, $[\zeta_{1,0}\to \zeta_{1,1}]$, and
$[\zeta_{1,1}\to \zeta_{2,0}]$.
The optional connection is $[\zeta_{1,1}\to \zeta_{3,0}]$
which is present in Figure~\ref{figStep3} (a).
\end{remark}

\section{Stability of ac-networks}
\label{sufcond}

In this section we derive sufficient conditions for asymptotic stability
of an ac-network $Y$, comprised
of equilibria $\xi_j$, $1\le j\le J$, and heteroclinic connections
$C_{ij}$. According to Lemma~\ref{lem1}, the connections
are either one or two-dimensional and any two-dimensional connection
is accompanied by two one-dimensional connections that
belong to a $\Delta$-clique.

In what follows we consider trajectories in the $\varepsilon$-neighbourhood
of $Y$ identifying this neighbourhood with the set $V$ from definition
\ref{as} of asymptotic stability. By $N_{\tilde\delta}(\xi_j)$ we denote
a family of $\tilde\delta$-neighbourhoods of $\xi_j$. We assume
$\tilde\delta$ to be sufficiently small (see subsection \ref{locglob}),
$\varepsilon\ll\tilde\delta$ and $\tilde\delta$ remains fixed as
$\varepsilon\to0$.

To study stability, the behaviour of a trajectory near a heteroclinic cycle or
network is approximated by local and global maps.
Namely, near an equilibrium
$\xi_j$ the local map $\bx_j^{out}=\phi_j\bx_j^{in}$ relates the point
$\bx_j^{out}$ where a trajectory $\Phi(\tau,x)$ exits the box
$N_{\tilde\delta}(\xi_j)$ with the point $\bx_j^{in}$ where it enters.
Outside the boxes, a connecting diffeomorphism $\psi_{ij}$ relates the point
$\bx_j^{in}$ where at $\tau=\tau_j^{in}$ a trajectory enters $N_{\tilde\delta}(\xi_j)$
with the point $\bx_i^{out}$ where at $\tau=\tau_i^{out}$ it exits
$N_{\tilde\delta}(\xi_i)$, assuming that for $\tau\in[\tau_i^{out},\tau_j^{in}]$ the
trajectory does not pass near any other equilibria.
We consider trajectories in the neighbourhoods $N_{\delta}(Y)$, for arbitrarily
small $\delta$ and use the following strategy:
we fix $\tilde\delta>0$ assuming $\tilde\delta$ to be small, but
$\delta\ll\tilde\delta$.

The results in Subsection~\ref{locglob} provide upper bounds for the distance
between the network and a trajectory exiting $N_{\tilde\delta}(\xi_j)$
as a function of the entry distance. The bounds involve exponents
$\rho_j$ that depend on the eigenvalues of $df(\xi_j)$. In Subsection~\ref{theo} we
prove that  trajectories remain close to the network when they enter
a subsequent neighbourhood and that they also remain close between neighbourhoods of consecutive equilibria.
We use these bounds to prove sufficient
conditions for asymptotic stability in the form of
 inequalities involving the
exponents. The results of Subsection~\ref{theo} are
applicable not only to ac-networks, but also in a more general setting,
e.g., for compact networks, if estimates for local maps are known.

To obtain estimates for ac-networks, instead of the distance $d(\cdot,\cdot)$
based on the maximum norm, it is convenient to use the complementary distance,
$\td(\cdot,\cdot)$, that we define below.

\begin{definition}\label{dist}
Consider a robust heteroclinic network $Y\subset\R^n$.
By definition of robust networks, for any connection
$C_{ij}\subset Y$ there exists a flow-invariant subspace $S_{ij}\supset C_{ij}$.
The {\em complementary distance} between a point $x\in\R^n$ and $Y$ is
$$
\td(x,Y)=\min d(x,S_{ij}),
$$
where the minimum is taken over all connections $C_{ij}$ in $Y$.
\end{definition}

\begin{remark}\label{dim2}
In general, a subspace $S_{ij}\supset C_{ij}$ is not unique and for different choices of
$S_{ij}$ one may possibly obtain different sufficient conditions for
asymptotic stability of a particular network. We choose $S_{ij}$ to be of
minimal dimension: for an ac-network the dimensions are $n_r+1$
(for one-dimensional connections) or $n_r+2$ (for two-dimensional ones),
where $n_r$ is the dimension of the radial eigenspace, the same for all
equilibria in a network.
\end{remark}

Note that $d(x,S_{ij})=\Pi_{ij}^\perp (x)$, where $\Pi_{ij}^\perp (x)$ denotes the orthogonal projection into $\R^n\ominus S_{ij}$, the orthogonal complement to $S_{ij}$ in $\R^n$. This is why we call it  {\em complementary distance}.

\subsection{Estimates for local   maps}
\label{locglob}

Choose $\tilde\delta$ sufficiently small so that
near $\xi_j$
the behaviour of trajectories can be approximated by
\begin{equation}\label{apploc}
\dot x_k=\lambda_{jk}x_k,
\end{equation}
where $\lambda_{jk}$ are eigenvalues of $\rd f(\xi_j)$ and $x_k$ are
coordinates in the local basis comprised of eigenvectors of $\rd f(\xi_j)$.
For the derivation of this approximation and its validity see, e.g.
\cite{Anosov,Samovol}.

For a heteroclinic cycle $X$ comprised of one-dimensional connections the
expression for a local map was derived in
\cite{KrupaMelbourne1,KrupaMelbourne2}; we give here a brief description.
Denote by $(u,v,w,\bz)$ the local coordinates in the basis comprised
of radial, contracting, expanding and transverse eigenvectors, respectively.
The incoming trajectory $\kappa_{ij}$ crosses the boundary of
$N_{\tilde\delta}(\xi_j)$ at a point $(u,v,w,\bz)=(u_0,v_0,0,0)$, the outgoing
trajectory $\kappa_{jk}$ crosses the boundary at $(u,v,w,\bz)=(u_1,0,w_1,0)$.
A trajectory close to the connection $\kappa_{ij}$ enters and exits
$N_{\tilde\delta}(\xi_j)$ at
\begin{equation}\label{inout}
\bx^{in}=(u_0+u^{in},v_0+v^{in},w^{in},\bz^{in})\hbox{ and }
\bx^{out}=(u_1+u^{out},v^{out},w_1+w^{out},\bz^{out}),
\end{equation}
respectively. Since $u^{in},v^{in},u^{out},w^{out}\ll u_0,v_0,u_1,w_1$, in the
expressions for $\bx^{in}$ and $\bx^{out}$ they are ignored. The radial
direction is ignored because it is irrelevant in the study of stability.
Therefore, the local map $(v^{out},\bz^{out})=\phi_{j}(w^{in},\bz^{in})$
can be approximated as
\begin{equation}\label{mapl0}
v^{out}=v_0w_1^{c/e}(w^{in})^{-c/e}\hbox{ and }
\bz^{out}=\bz^{in}w_1^{\bt/e}(w^{in})^{-\bt/e},
\end{equation}
where $c,e,\bt$ are the contracting, expanding and transverse eigenvalues,
respectively.

Since the radial variable in
$\bx^{out}$ depends on $w^{in}$ as $u^{out}=u_0w_1^{c/e}(w^{in})^{-r/e}$,
we use the robust distances as the measure of closeness of $\Phi(\tau,x)$ to
a cycle, or to a network.
Let $\td^{in}=\td(\bx^{in},X)=|(w^{in},\bz^{in})|$ and
$\td^{out}=\td(\bx^{out},X)=|(v^{out},\bz^{out})|$ be the robust distances between
the entry and exit points of $\Phi(\tau,x)$ and $X$.
From (\ref{mapl0}) we obtain the estimate
\begin{equation}\label{estl0}
\td^{out}<A_j(\td^{in})^{\rho_j},\hbox{ where }
\rho_j=\min(-c/e,1+\min_s(-t_s/e))
\end{equation}
and $A_j>0$ is a constant independent of $w^{in}$ and $\bz^{in}$.

\begin{remark}\label{rem41}
Let $\xi_j\in Y$, where $Y$ is an ac-network, be an equilibrium with
one-dimensional incoming and outgoing connections. Then
the approximations (\ref{mapl0}) and the estimate
(\ref{estl0}) for the local map near $\xi_j$ hold true. Note
that for an ac-network all transverse eigenvalues are negative.
\end{remark}

From now on we deal with the case when there are 2-dimensional connections.

Recall that in a $\Delta$-clique (Figure~\ref{fig:Delta-clique}) a b-point (beginning) is connected to both an m-point (middle) and an e-point (end) with an additional connection from the  m-point to the e-point.

\begin{figure}[h]
\begin{center}
\parbox{40mm}{
\begin{center}
\includegraphics[width=39mm]{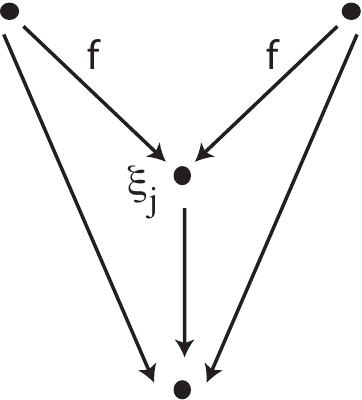}\\
Lemma~\ref{lem43}
\end{center}
}\qquad
\parbox{40mm}{
\begin{center}
\includegraphics[width=39mm]{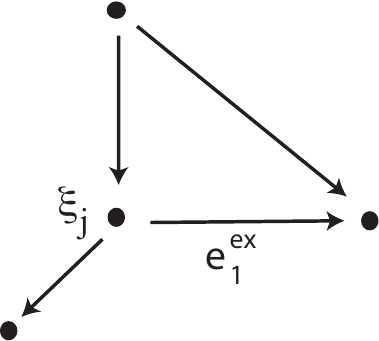}\\
Lemma~\ref{lem44}
\end{center}
}\\
\parbox{100mm}{
\begin{center}
\includegraphics[width=45mm]{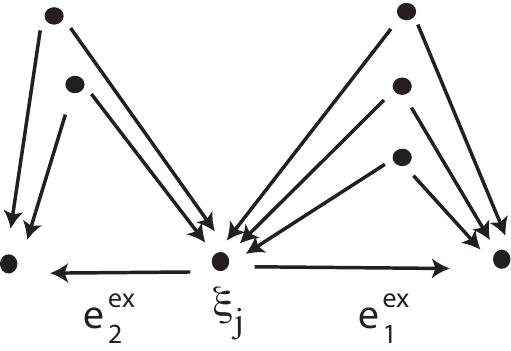}
\qquad
\includegraphics[width=45mm]{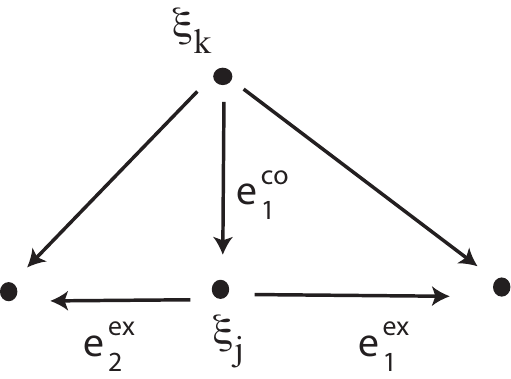}\\
Lemma~\ref{lem45}
\end{center}
}

\end{center}
\caption{The relevant structure in the networks of
Lemmas~\ref{lem43}, \ref{lem44} and \ref{lem45}.  The  right-hand figure for
Lemma~\ref{lem45} shows a non ac-network.}\label{fig:sec41}
\end{figure}

We divide the estimates into lemmas taking into account whether the equilibrium
is, or is not, an m-point for a $\Delta$-clique and on whether it has one or two
expanding eigenvectors.
 Estimates around a node that is a b-point and/or an
e-point appear in Lemma~\ref{lem41}.
The cases of an m-point with only one expanding eigenvalue are covered
in Lemmas~\ref{lem42} and \ref{lem43} and the ones with two expanding
eigenvalues in Lemmas~\ref{lem44} and \ref{lem45}
(see Figure~\ref{fig:sec41}). The combinations cover all possibilities for
a node in an ac-network.

\bigbreak

Let $\xi_j$ be an equilibrium in an ac-network $Y$. The
notions of radial, contracting,
expanding and transverse eigenvalues of the Jacobian matrix at an equilibrium are given in Subsection~\ref{sec:robust_network}.
Denote by $c_s$, $1 \leq s \leq S$, the contracting eigenvalues of $\rd f(\xi_j)$;  by $e_p$, $1\le p \le P$, the expanding eigenvalues; finally, by $t_q$, $1 \le q \le Q$, the transverse eigenvalues.
Note that $P$ is either 1 or 2. We write
$$
c=\max_{1\le s\le S}c_s \;\;\; \mbox {  and   } \;\;\; e=\max_{1\le p\le P}e_p.
$$
We use $\be^{co}_s$ to denote the contracting eigenvectors ($1\le s \le S$) and $\be^{ex}_p$ ($p= 1$ or $2$) the expanding eigenvectors.
We use the superscripts $.^{in}$ and $.^{out}$ to indicate the quantities at the moments when a trajectory $\Phi(\tau,x)$
enters and exits $N_{\tilde\delta}(\xi_j)$, respectively. As above, the corresponding time moments are $\tau^{in}$ and $\tau^{out}$ and
we define $\td^{in}=\td(\bx^{in},Y)$ and $\td^{out}=\td(\bx^{out},Y)$.

\begin{lemma}\label{lem41}
Suppose that an equilibrium $\xi_j\in Y$ is not an m-point for
any of the $\Delta$-cliques of $Y$. Let $\xi_j$ have $S$ contracting eigenvectors
and $P$ expanding
eigenvectors. Then
\begin{equation}\label{estl}
\td^{out}<A_j(\td^{in})^{\rho_j},\hbox{ where }
\rho_j=\min(-c/e,1+\min_q(-t_q/e)),\
A_j \in \R.
\end{equation}
\end{lemma}

\proof
Since $\xi_j$ is not an m-point, the expanding eigenvectors $\be^{ex}_p$
belong to the orthogonal complement to any of $P_d^{in}$, where
$P_d^{in}$ are the invariant subspaces that the incoming connections
belong to. Therefore, $\max_{1\le p\le P}|w_p^{in}|<\td^{in}$. If $P=2$ then the
outgoing connections (short and f-long connections of a $\Delta$-clique)
intersect with $N_{\tilde\delta}(\xi_j)$ along the faces $w_1=\tilde\delta$ and
$w_2=\tilde\delta$. For the trajectory that exits through the face
$w_1=\tilde\delta$ the time of flight $\tau^*=\tau^{out}-\tau^{in}$ satisfies
$$\re^{\tau^*}=\tilde\delta^{1/e_1}|w_1^{in}|^{-1/e_1}.$$
For the trajectory that exits through the face $w_2=\tilde\delta$ it satisfies
$$\re^{\tau^*}=\tilde\delta^{1/e_2}|w_2^{in}|^{-1/e_2}.$$
Therefore, we have
$$\re^{\tau^*}>C_1(\td^{in})^{-1/e},$$
where $C_1$ is a constant. If $P=1$ then the above estimate evidently holds true.

There exists a constant $C_2>0$ such that $v_{s,0}<C_2$ for all
incoming connections. Therefore, for the contracting coordinates
\begin{equation}\label{esl41}
v^{out}_s=v_{s,0}\re^{c\tau^*}<C_2(C_1)^c(\td^{in})^{-c/e}.
\end{equation}
For the transverse ones we have
\begin{equation}\label{esl42}
z^{out}_q=z_{q,0}\re^{t_q\tau^*}<(C_1)^{t_q}(\td^{in})^{1-t_q/e}.
\end{equation}
The statement of the lemma follows from (\ref{esl41}), (\ref{esl42}) and
the definition of the robust distance.
\qed

\begin{lemma}\label{lem42}
Suppose that an equilibrium $\xi_j\in Y$ is m-point for only one
$\Delta$-clique of $Y$, $\xi_j$ has one contracting eigenvector
and  only one expanding, which are the f-long and s-long vectors
of the $\Delta$-clique. Then
the trajectories passing through $N_{\tilde\delta}(\xi_j)$ satisfy
\begin{equation}\label{estl1}
\td^{out}\le \td^{in}.
\end{equation}
\end{lemma}

\proof
The expanding and contracting eigenvectors of $\xi_j$ belong to the invariant
subspace $P_j$ that contains the $\Delta$-clique. Therefore, $P_j^{\perp}$
is spanned by the transverse eigenvectors of $\rd f(\xi_j)$. Since the
transverse coordinates satisfy $z^{out}_q=z_{q,0}\re^{t_q\tau^*}$, where again $\tau^*=\tau^{out}-\tau^{in}$, and
all the transverse eigenvalues $t_q$ are negative the statement of the lemma
holds true.
\qed

\begin{lemma}\label{lem43}
Suppose that an equilibrium $\xi_j\in Y$ is m-point for several
(one or more) $\Delta$-cliques of $Y$ and it has one expanding eigenvector.
For $1\le s\le S_0$ let $\be^{co}_s$ be the f-long vector for a $\Delta$-clique, while the remaining contracting eigenvectors are not f-long.
Then
\begin{equation}\label{estl}
\td^{out}<A_j(\td^{in})^{\rho_j},\hbox{ where }
\rho_j=\min(-c/e,1)\hbox{ and } A_j \in \R.
\end{equation}
\end{lemma}

\proof
First we consider trajectories that enter the $\delta$-neighbourhood
of a $\Delta$-clique.
Since all eigenvalues of $\rd f(\xi_j)$, except for one expanding, are negative, then
a trajectory $\Phi(\tau,x)$ that at $\tau=\tau^{in}$ belongs to the $\delta$-neighbourhood
of a $\Delta$-clique remains in this neighbourhood
for $\tau\in[\tau^{in},\tau^{out}]$ and satisfies
\begin{equation}\label{esl46}
\td^{out}<\td^{in}.
\end{equation}

If a trajectory does not belong to the
$\delta$-neighbourhood of any of $\Delta$-cliques then it satisfies
$|w^{in}|<\td^{in}$. Hence, the time of flight $\tau^*=\tau^{out}-\tau^{in}$ satisfies
$$\re^{\tau^*}>C_1(\td^{in})^{-1/e}.$$
(This can be shown by the same arguments as the ones used in the proof of Lemma
\ref{lem41}.) Let $C_2>0$ be a constant such that $v_{s,0}<C_2$ for all
incoming connections. Then
\begin{equation}\label{esl44}
v^{out}_s=v_{s,0}\re^{c\tau^*}<C_2(C_1)^c(\td^{in})^{-c/e},
\end{equation}
where $c=\max_{S_0< s\le S}c_s$.
For the transverse coordinates we have
\begin{equation}\label{esl45}
z^{out}_q=z_{q,0}\re^{t_q\tau^*}<(C_1)^{t_q}(\td^{in})^{1-t_q/e}.
\end{equation}
The statement of the lemma follows from (\ref{esl46})-(\ref{esl44}) and
the definition of the robust distance.
\qed

The following lemma provides bounds for quantities used in the proof of Lemma~\ref{lem44}.

\begin{lemma}\label{lem440}
Suppose that $\lambda_1<0$ and $\lambda_2,q,\varepsilon>0$. Then
$$\min(q^{\lambda_1},\varepsilon q^{\lambda_2})\le
\varepsilon^{\lambda_1/(\lambda_1-\lambda_2)}.$$
\end{lemma}

\proof
$$
\renewcommand{\arraystretch}{1.2}
\begin{array}{l}
\hbox{If $q^{\lambda_1}\le\varepsilon q^{\lambda_2}$ then
$q^{\lambda_1-\lambda_2}\le\varepsilon,$ which implies that
$q^{\lambda_1}\le\varepsilon^{\lambda_1/(\lambda_1-\lambda_2)}$.}\\
\hbox{If $\varepsilon q^{\lambda_2}\le q^{\lambda_1}$ then
$q^{\lambda_2-\lambda_1}\le\varepsilon^{-1}$, which implies that
$\varepsilon q^{\lambda_2}\le\varepsilon^{{\lambda_1}/(\lambda_1-\lambda_2)}$.}
\end{array}
$$
\qed

\begin{lemma}\label{lem44}
Suppose that an equilibrium $\xi_j\in Y$ has one contracting eigenvector,
$\be^{co}$, two expanding eigenvectors, $\be^{ex}_1$ and $\be^{ex}_2$,
$\xi_j$ is an m-point for just one $\Delta$-clique,
$\be^{co}$ is  the f-long vector of the $\Delta$-clique and
$\be^{ex}_1$ is the s-long vector of the $\Delta$-clique.
Then
\begin{equation}\label{estle44}
\td^{out}<A_j(\td^{in})^{\rho_j},\hbox{ where }
\rho_j={c\over c-e_2}, \  A_j \in \R.
\end{equation}
\end{lemma}

\begin{figure}[h]
\centerline{
\includegraphics[width=9cm]{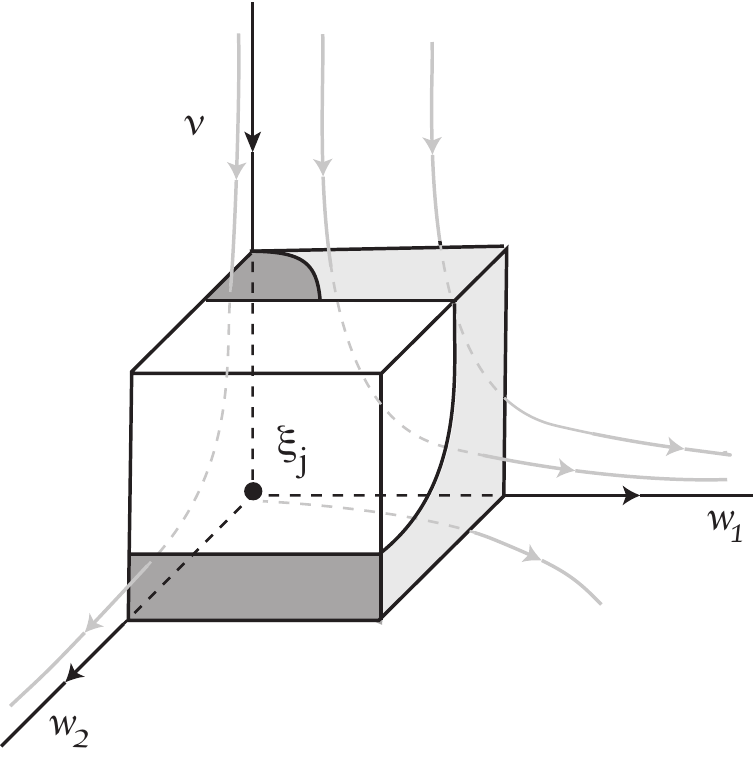}
}
\caption{The neighbourhood of $\xi_j$ considered in Lemma \ref{lem44}
 in the case where the largest expanding eigenvalue is in the $w_1$ direction.
Trajectories in $Y$ are black, while the ones near $Y$ are grey.
The dark grey regions indicate the entry and exit points of the trajectories that
exit through the face $w_2=\tilde\delta$ while the light grey indicates the ones
that exit through the face $w_1=\tilde\delta$.}

\end{figure}

\proof
Trajectories that enter $N_{\tilde\delta}(\xi_j)$ belong to the
$\delta$-neighbourhood of the $\Delta$-clique, hence initially
$w^{in}_2<\td^{in}$. Therefore, the trajectories that exit through the face
$w_2=\tilde\delta$ satisfy
\begin{equation}\label{est01}
v^{out}<C(\td^{in})^{-c/e_2}\hbox{ and }
z^{out}_q<C(\td^{in})^{1-t_q/e_2}.
\end{equation}
(This follows from  arguments similar to the ones employed in the proofs
of Lemmas \ref{lem41} and \ref{lem43}.)

For trajectories that exit through the face $w_1=\tilde\delta$ the time
of flight satisfies
$$w^{in}_1\re^{\tau^*}=\tilde\delta.$$
At the exit point
$$v^{out}=v_{0}\re^{c\tau^*}\hbox{ and }w^{out}_2=w_2^{in}\re^{e_2\tau^*}.$$
Therefore (see Lemma \ref{lem440})
\begin{equation}\label{est02}
\min(v^{out},w^{out}_2)<C(\td^{in})^{c/(c-e_2)}.
\end{equation}
The statement of the lemma follows from (\ref{est01}), (\ref{est02}) and
the fact that $-c/e_2>c/(c-e_2)$ and that
$$z^{out}_q=z^{in}_q\re^{t_q\tau^*}\hbox{ where }t_q<0.$$
\qed

\begin{lemma}\label{lem45}
Suppose that an equilibrium $\xi_j\in Y$ has several contracting eigenvectors and two expanding eigenvectors.
Let the eigenvectors $\be^{co}_1,...,\be^{co}_q$
be f-long vectors of $\Delta$-cliques for which $\be^{ex}_1$ is the s-long vector.
Assume also that the eigenvectors $\be^{co}_h,...,\be^{co}_d$, where $h\le q+1$,
 are f-long vectors of $\Delta$-cliques where $\be^{ex}_2$ is the s-long vector. Then
\begin{equation}\label{estle45}
\renewcommand{\arraystretch}{1.2}
\begin{array}{l}
\td^{out}<A_j(\td^{in})^{\rho_j},\hbox{ where }
\rho_j=\min(a_1,a_2,a_3),\
a_1=-\max_{d+1\le s\le S} c_s/\max_{p=1,2} e_p\\
a_2=\tilde c_1/(\tilde c_1-e_1),\ a_3=\tilde c_2/(\tilde c_2-e_2),\
\tilde c_1=\max_{1\le s\le q} c_s,\
\tilde c_2=\max_{h\le s\le d} c_s, \ A_j \in \R.
\end{array}
\end{equation}
\end{lemma}

\proof
Proceeding as in the proof of Lemma \ref{lem44}, we first obtain estimates of
$\td^{out}$ for trajectories that when they enter $N_{\tilde\delta}(\xi_j)$
are not in the $\delta$-neighbourhoods of any $\Delta$-cliques, then we consider
trajectories in each of the $\Delta$-cliques individually.
\qed

\begin{remark}\label{acnetw}
Existence of a common f-long connection for two $\Delta$-cliques means that  the
$\Delta$-cliques have a common b-point whose unstable manifold has dimension three or larger,
as in the last panel of Figure~\ref{fig:sec41}.
Since the dimension of the unstable manifold of an equilibrium in an ac-network
is one or two, this implies that for an ac-network we have $h=q+1$.
We prove the lemma in its more general form so that it also applies
to other heteroclinic networks.
\end{remark}

\subsection{Sufficient conditions for asymptotic stability}\label{theo}

In this section we prove a theorem providing sufficient conditions for
asymptotic stability of certain heteroclinic networks. The networks
that are considered in theorems\footnote
{We state two theorems providing sufficient conditions for asymptotic
stability. Since the proofs are similar, we prove only one of them.}
and lemmas of this section are robust heteroclinic networks in
(\ref{eq:ode}), comprised of a finite number of hyperbolic equlibria $\xi_j$
and heteroclinic connections $C_{ij}$ that belong to flow-invariant
subspaces $S_{ij}$. In each lemma and theorem we explicitly state the
other assumptions that are made. The ac-networks satisfy these assumptions,
therefore the lemmas and theorems are applicable to them.

In this subsection we consider a set $N_{\tilde\delta}(\xi_j)$ bounded by $|a_m|=\tilde\delta$,
where the $a_m$ coordinate is taken in a basis comprised of eigenvectors of $df(\xi_j)$.
Thus, $N_{\tilde\delta}(\xi_j)$  is not a box, as it is for ac-networks, but a
bounded region with boundary denoted $\partial N_{\tilde\delta}(\xi_j)$.
 Recall that we consider trajectories in
a $\delta$-neighbourhood of the network, where $0<\delta\ll\tilde\delta$ and
$\tilde\delta$ is fixed and small. The times $\tau^{in}$ and $\tau^{out}$ are
taken with reference to $N_{\tilde\delta}(\xi_j)$.
Lemmas~\ref{lem51} and \ref{lem51nn} are concerned
with bounds for the distance from a trajectory to the network for times in
$[\tau^{in},\tau^{out}]$, that is, when the trajectory is near a node.
By Lemma~\ref{lem51} smallness of the distance $d(\Phi(\tau),Y)$
at $\tau=\tau^{in}$ implies its smallness for any
$\tau\in[\tau^{in},\tau^{out}]$, while by Lemma \ref{lem51nn} the smallness
at some $\tau\in[\tau^{in},\tau^{out}]$ implies its smallness at
$\tau=\tau^{out}$. In Lemmas~ \ref{lem52} and \ref{lem52n} similar results
are proven along a connection from one node to the next one, i.e. for
the interval $[\tau^{out},\tau^{in}]$. According to Lemma~\ref{lem51n}
in the study of stability one can use the robust distance $\td(y,Y)$
instead of $d(y,Y)$. In Theorems \ref{th_as} we prove sufficient
conditions for asymptotic stability that are inequalities involving exponents
$\rho_j$ of the local maps obtained in Section~\ref{locglob}. Theorem
\ref{th_as2} concerns more cumbersome conditions of asymptotic stability
that involve slightly different exponents for estimates of local maps.

\begin{lemma}\label{lem51}
Given a clean robust heteroclinic network $Y$, there exist $s>0$ and
$\varepsilon'>0$ such that for any $0<\varepsilon<\varepsilon'$
$$d(\Phi(\tau^{in},x),Y)<\varepsilon^s\hbox{ implies that }
d(\Phi(\tau,x),Y)<\varepsilon\hbox{ for all }\tau\in[\tau^{in},\tau^{out}].
$$
\end{lemma}

\proof
Since near $\xi_j$ a flow of (\ref{eq:ode}) can be approximated by
a linearised system, in $N_{\tilde\delta}(\xi_j)$ the difference
$\Delta\bx=\bx(\tau,\bx_0)-\by(\tau,\by_0)$ between two trajectories is
\begin{equation}\label{eq420}
\Delta\bx=(\bx_0-\by_0)\re^{\Lambda \tau},
\end{equation}
where $\Lambda=(\lambda_1,...,\lambda_n)$ are the eigenvalues of $\rd f(\xi_j)$.
For simplicity we assume that all eigenvalues are real. In the case of complex
eigenvalues the proof is similar and in the expressions for $s_j$
the eigenvalues should be replaced by their real parts.

Let $\bx_0$ be the point where $\Phi(\tau,x)$ enters $N_{\tilde\delta}$ and
$\by_0$ be the nearest to $\bx_0$ point in $Y$. Split contracting eigenvectors
of $\rd f(\xi_j)$ into two groups: belonging to the subspace
$S_{ij}\supset C_{ij}\ni\by_0$, and those that do not belong.
The respective coordinates are denoted by $\bv$ and $\bq$ and the associated
eigenvalues by $\bc_v$ and $\bc_q$. The
expanding eigenvectors are splitted similarly, into belonging to $S_{ij}$
and do not. The respective coordinates are denoted by $\bh$ and $\bw$,
associated eigenvalues by $\be_h$ and $\be_w$.
Due to (\ref{eq420})
$$d(\Phi(\tau,\bx_0),\Phi(\tau,\by_0))=
\max(|\Delta\bu\re^{\br \tau}|,|\Delta\bv\re^{\bc_v \tau}|,|\Delta\bh\re^{\be_h \tau}|,
|\bq_0\re^{\bc_q \tau}|,|\bw_0\re^{\be_w \tau}|,|\bz_0\re^{\bt \tau}|),$$
where we have denoted
$$\bx_0=(\bu_0,\bv_0,\bh_0,\bq_0,\bw_0,\bz_0),$$
$$\by_0=(\bu_0-\Delta\bu,\bv_0-\Delta\bv,\bh_0-\Delta\bh,0,0,0).$$

Similarly we write
$$d(\Phi(\tau,\bx_0),W^u(\xi_j))=
\max(|\bu_0\re^{\br \tau}|,|\bv_0\re^{\bc_v \tau}|,
|\bq_0\re^{\bc_q \tau}|,|\bz_0\re^{\bt \tau}|).$$
Since all transverse eigenvalues are negative, using Lemma (\ref{lem440})
we obtain that if
$$d(\bx_0,\by_0)=\min(
|\Delta\bu|,|\Delta\bv|,|\Delta\bh|,|\bq_0|,|\bh_0|,|\bz_0|)<\varepsilon^s$$
then for all $\tau\in[\tau^{in},\tau^{out}]$
$$d(\Phi(\tau,\bx_0),Y)<
\min\biggl(d(\Phi(\tau,\bx_0),\Phi(\tau,\by_0)),d(\Phi(\tau,\bx_0),W^u(\xi_j)\biggr)
<\varepsilon^{sd_j},$$
where $d_j=c'/(c'-e)$, $c'=\max(r,c)$, $r=\max_{1\le k\le K}r_k$ and
$r_k$ are the radial eigenvalues of $\xi_j$. (Recall that we do not assume
that $Y$ is an ac-network, hence its equilibria can have several radial
eigenvalues). Taking $s$ to be the maximum of $1/d_j$ we prove the lemma.
\qed

\begin{lemma}\label{lem51nn}
Let $Y$ be the network from Lemma \ref{lem51}.
Then there exist $s>0$ and $\varepsilon'>0$ such that for any
$0<\varepsilon<\varepsilon'$
$$d(\Phi(\tau,x),Y)<\varepsilon^s\hbox{ for some }\tau\in[\tau^{in},\tau^{out}]
\hbox{ implies that }d(\Phi(\tau^{out},x),Y)<\varepsilon.
$$
\end{lemma}

The proof is similar to the proof of Lemma \ref{lem51} and we do not present it.

In the following lemma $Y$ is a heteroclinic network comprised of
equilibria $\xi_j$, $1\le j\le J$, and connecting trajectories, and
$$
\cN=\{N_{\tilde\delta}(\xi_1),...,N_{\tilde\delta}(\xi_J)\}
$$
is a family of disjoint neighbourhoods of the equilibria.
Recall that $\Phi(\tau,x)$ is a trajectory through $x$ such that $\Phi(0,x)=x$. 
Let $\Phi((\tau_1,\tau_2),x)$  be the part of the trajectory of $x$ between times $\tau_1$ and $\tau_2$.
For a given $x\not\in N_{\tilde\delta}(\xi_{j})$,
for any $1\le j\le J$, we denote by $\cF(x,\cN)$ the first entrance point
of the trajectory $\Phi(\tau,x)$ into one of $N_{\tilde\delta}(\xi_k)$, that is,
\begin{itemize}
\item[]  $\cF(x,\cN)=\Phi(\tau^{in},x)\in  \partial N_{\tilde\delta}(\xi_k)$,
 \item[]  $\Phi((\tau^{in},\tau^{in}+q),x)  \in  N_{\tilde\delta}(\xi_k)
  \hbox{ for some }q>0$, and
\item[]  $\Phi(\tau,x)  \notin  N_{\tilde\delta}(\xi_{j})\hbox{  for any }
\tau\in(0,\tau^{in})\hbox{  and }1\le j\le J$.
\end{itemize}
We denote the time of the first entry by $\Theta(x,\cN)=\tau^{in}$.
For a trajectory that never enters any neighbourhood in $\cN$ for all positive $\tau$ we write $\Theta(x,\cN)=\infty$. Note that for a point $x \in  \partial N_{\tilde\delta}(\xi_j)$ we may have $\cF(x,\cN)=x$ and $\Theta(x,\cN)=0$ if the trajectory through $x$ enters $N_{\tilde\delta}(\xi_j)$.

\begin{lemma}[Lemma 4.13]\label{lem52}
Given a compact
 heteroclinic network $Y$
and a family of disjoint neighbourhoods $\cN$
of its equilibria $\xi_j$, $1\le j\le J$, then
for any equilibrium $\xi_*\in Y$
there exists $B>0$ and $\varepsilon_0>0$
such that for any $0<\varepsilon<\varepsilon_0$
$$x\in\partial N_{\tilde\delta}(\xi_*) \hbox{ and }
d(x,Y)<\varepsilon\hbox{ imply that }$$
$$\Theta(x,\cN)\neq\infty\hbox{ and }
d(\Phi(\tau,x),Y)<B\varepsilon\hbox{ for all }\tau\in[0,\Theta(x,\cN)].
$$
\end{lemma}

\proof
The intersection $Y\cap \partial N_{\tilde\delta}(\xi_*)$ can be decomposed as
$$
Y\cap \partial N_{\tilde\delta}(\xi_*)=[Y\cap \partial N_{\tilde\delta}(\xi_*)\cap W^u(\xi_*)]\cup
[Y\cap \partial N_{\tilde\delta}(\xi_*)\cap W^s(\xi_*)].
$$
For an initial condition $x$ that is close to the latter set in
the union we have $\Theta(x,\cN)=0$, therefore the statement of the
lemma is trivially satisfied with $B=1$ and any $0<\varepsilon_0<\tilde\delta$.

In the remainder of this proof we consider only
the initial conditions $x$ that are close to
$\cW = Y\cap \partial N_{\tilde\delta}(\xi_*)\cap W^u(\xi_*)$. 
The set $\cW$ is compact due to the compactness of $Y$ and
$\partial N_{\tilde\delta}(\xi_*)$.
Evidently for $x$ close to $\cW$ we have $d(x,Y)=d(x,Y\cap W^u(\xi_*))$.

For $x\in\partial N_{\tilde\delta}(\xi_*)$ such that
$d(x,Y)<\varepsilon$, denote by $y$ and $y'$ the 
points in $Y\cap\partial N_{\tilde\delta}(\xi_*)$ and $Y$, respectively,
nearest to $x$.
 The points $y$ and $y'$ may be different since $W^u(\xi_*)$ is not flat and the faces of the box $N_{\tilde\delta}(\xi_*)$ may not  be orthogonal to $W^u(\xi_*)$.
We have $|x-y'|=|x-y|/\cos(\theta)$, where $\theta$ is the angle
between the vectors $(x-y)$ and $(x-y')$.
Let $\phi$ be the smallest angle between a face of
$\partial N_{\tilde\delta}(\xi_*)$ and the unstable subspace of $\rd f(\xi_*)$.
\footnote{For an ac-network,  $\phi=\pi/2$.}
Since $W^u(\xi_*)$ is tangent to the unstable subspace of $\rd f(\xi_*)$, for $\tilde\delta\to0$ the
angle between the vector $(y-\xi_*)$ and the unstable subspace vanishes,
implying that in this case $\theta\le \pi/2-\phi$.
Therefore, for sufficiently small $\tilde\delta$ the condition
$d(x,Y)<\varepsilon<\tilde\delta$ implies existence
of $y\in \cW$
such that 
$|x-y|<C\varepsilon$, 
where $C=2/\sin\phi$.
 Below we assume that
$\tilde\delta$ is such sufficiently small.

Introduce another family of smaller neighbourhoods
$$
\cN'=\{N_{\tilde\delta/2}(\xi_1),...,N_{\tilde\delta/2}(\xi_J)\}.
$$
For  the set of equilibria  $\xi_s$, $1\le s\le S$,  in $Y$ for which there is a connection $[\xi_*\to\xi_{s}]$,
decompose $\cW$ as
$$
\cW=\cup_{1\le {s}\le S}\cW_{s}\quad\hbox{ where }\quad
\cW_{s}=\{y\in\cW~:~\cF(y,\cN')\in
\partial N_{\tilde\delta/2}(\xi_{s})\}.
$$
Due to the smoothness of $f$, for small $|x-y|$
the difference $\Phi(\tau,x)-\Phi(\tau,y)$ satisfies
$$
|\Phi(\tau,x)-\Phi(\tau,y)|<B(\tau,y)|x-y|,
$$
for some bound $B$ depending only on the time $\tau$ and
the initial condition $y$. 

Since $\overline{\cW}_s$ is compact, there exist
$\varepsilon_s>0$ and $0<B_s<\infty$ such that
for all
$y\in\overline{\cW}_s$ and $ x\in\partial N_{\tilde\delta}(\xi_*)$ satisfying
$|x-y|<\varepsilon_s$ and $0<\tau<\Theta(y,\cN')$
we have
\begin{equation}
|\Phi(\tau,x)-\Phi(\tau,y)|<B_s|x-y|.
\label{eq1}
\end{equation}

Let
\begin{equation}\label{eq2}
\varepsilon_0=\frac{1}{C}\min_{1\le s \le S}
(\min(\frac{\tilde\delta}{2B_s},
\varepsilon_s))\quad\hbox{ and }\quad
B=\max_{1\le s \le S}B_s.
\end{equation}
Below we prove that a trajectory $\Phi(\tau,x)$, where $x\in\partial N_{\tilde\delta}(\xi_*)$
and $|x-y|<C\varepsilon_0$ for some $y\in\cW_s$, satisfies
$$
\Phi(\Theta(y,\cN'),x)\in N_{\tilde\delta}(\xi_s).
$$
In fact, we have
$$
|\Phi(\Theta(y,\cN'),x)-\Phi(\Theta(y,\cN'),y)|<B_s |x-y|<
B_sC \varepsilon_0\le B_sC \frac{\tilde{\delta}}{2B_sC}
=\frac{1}{2}\tilde{\delta}.
$$
Since 
$\Phi(\Theta(y,\cN'),y)\in\partial{N_{\tilde\delta/2}(\xi_s)}$
then $|\Phi(\Theta(y,\cN'),y)-\xi_s | \le \tilde\delta/2$ and we get
$$
|\Phi(\Theta(y,\cN'),x)-\xi_s | \le |\Phi(\Theta(y,\cN'),x)-\Phi(\Theta(y,\cN'),y)| + |\Phi(\Theta(y,\cN'),y)-\xi_s |
<  \frac{1}{2}\tilde{\delta}+\frac{1}{2}\tilde{\delta}.
$$
Hence, $\Phi(\Theta(y,\cN'),x) \in N_{\tilde\delta}(\xi_s)$ and
therefore $\Theta(x,\cN)<\Theta(y,\cN')<\infty$.

Due to (\ref{eq1}) and (\ref{eq2}), for any
$0<\tau<\Theta(x,\cN)<\Theta(y,\cN')$ we have
$$
d(\Phi(\tau,x),Y) \le |\Phi(\tau,x)-\Phi(\tau,y)|<B_s |x-y|<B|x-y|,
$$
finishing the proof. \qed

\begin{lemma}\label{lem52n}
Let $Y$ be the network from Lemma \ref{lem52}.
For any heteroclinic connection $C_{ij}$ there exists $B_{ij}>0$ such that
$$d(\Phi(\tau,x),Y)<\varepsilon\hbox{ for some }\tau\in[\tau^{out},\tau^{in}]\hbox{ implies that }
d(\Phi(\tau^{in},x),Y)<B_{ij}\varepsilon.$$
\end{lemma}

The proof is similar to the proof of Lemma \ref{lem52} and we do not present it.

\begin{lemma}\label{lem51n}
Let $Y$ be the network from Lemma \ref{lem51} and $\xi_j\in Y$ be such that
$\bar W^u(\xi_j)\cap W^u(\xi_i)=\xi_j$ or $\varnothing$ for any $i\ne j$.
Then there exist $s>0$ and
$\varepsilon'>0$ such that for any $0<\varepsilon<\varepsilon'$
$$\td(\Phi(\tau^{in},\bx),Y)<\varepsilon^s\hbox{ implies that }
|\bu^{out}|<\varepsilon.$$
\end{lemma}

\proof
The conditions of the lemma ensure the absence of $\bh$ coordinated introduced
in the proof of Lemma \ref{lem51}. Since $|\bw|<\varepsilon^s$, for a trajectory
exiting through the face $w_p=\tilde\delta$ we have
$$|\bu^{out}|<|\bu^{in}|\re^{-r\tau^*}<\tilde\delta^{1+r/c_p}
\varepsilon^{-rs/c}.$$
Therefore, for $s=-c/2r$ the statement of the lemma holds true.
\qed

\begin{remark}\label{rem51}
In an ac-network any equilibrium which is not an m-point for
a $\Delta$-clique satisfies the conditions of Lemma \ref{lem51n}.
\end{remark}

\begin{theorem}\label{th_as}
Let $Y$ be a
clean heteroclinic network such that for any equilibrium
the nearby trajectories satisfy
$$
d(\Phi(\tau^{out},x),Y)<A_j(d(\Phi(\tau^{in},x),Y))^{\rho_j}
\hbox{ for some }A_j>0.
$$
Suppose also that for any heteroclinic cycle $X\subset Y$,
 involving $\xi_1,\ldots,\xi_m$, we have
$$\rho(X)>1
\qquad\mbox{where}\qquad
\rho(X)=\prod_{1\le j\le m}\rho_j.
$$
Then $Y$ is asymptotically stable.
\end{theorem}

\proof
According to the definition of asymptotic stability, to prove that $Y$ is
a.s., for any given $\delta>0$ we should find $\varepsilon>0$ such that
any trajectory $\Phi(\tau,x)$, with $d(x,Y)<\varepsilon$,  satisfies
\begin{itemize}
\item [(i)] $d(\Phi(\tau,x),Y)<\delta$ for any $\tau>0$;
\item [(ii)] $\lim_{\tau\to\infty}d(\Phi(\tau,x),Y)=0$.
\end{itemize}
Since $Y$ is clean, by Lemmas \ref{lem51}-\ref{lem52n}
it is sufficient to check (i) only at time
instances when $\Phi(\tau,x)$ enters the boxes $N_{\tilde\delta}(\xi_j)$.
For a trajectory which is staying in a $\delta$-neighbourhood
of $Y$
for all $\tau>0$ there are two possibilities: it can either be attracted to an
equilibrium $\xi_j\in Y$, or it can make infinite number of crossings of the
boxes $N_{\tilde\delta}(\xi_j)$. Therefore, the condition (ii) also may be
checked only at the time instances when the trajectory enters the boxes.
So, we consider only times when $\Phi(\tau,x)$ enters a box
and approximate a trajectory $\Phi(\tau,x)$ by a sequence of maps
$g_{ij}=\psi_{ij}\phi_i(x)$, where $x\in\partial N_{\tilde\delta}(\xi_i)$.
For any connection $C_{ij}$ Lemma \ref{lem52} ensures the existence of
$B_{ij}>0$ such that
$$d(\Phi(\tau^{in},x),Y)<B_{ij}d(\Phi(\tau^{out},x),Y).$$

Since for any  cycle $X=[\xi_1\to...\to\xi_m\to\xi_1]$
we have $\rho(X)>1$, there exists $q>0$ such that
$$\tilde\rho(X)=\prod_{1\le j\le m}\tilde\rho_j>1,
\hbox{ where }\tilde\rho_j=\rho_j-q.$$

Given a path $\cP=[\xi_1\to...\to\xi_m]$ define
$$
\rho(\cP)=\prod_{1\le j\le m}\rho_j \;\;\; \mbox{ and  } \;\;\;  \tilde\rho(\cP)=\prod_{1\le j\le m}\tilde\rho_j
$$
Denote $A=\max_jA_j$ and $B=\max_{ij}B_{ij}$, where the maxima are
taken over all equilibria and all connections in $Y$.
We choose $\varepsilon$ according to the following rules:
\begin{itemize}
\item [I)] $0<\varepsilon<\delta$;
\item [II)] $\varepsilon^{(\tilde\rho(X)-1)\tilde\rho(\cP)}<1/2$ for any
heteroclinic cycle $X\subset Y$ and any path $\cP\subset Y$;
\item [III)] $\varepsilon^{\tilde\rho(\cP)}<\tilde\varepsilon$
for any path ${\cal P}\subset Y$, where $0<\tilde\varepsilon<\delta$ satisfies
$\tilde\varepsilon^qAB<1$.
\end{itemize}

Consider a trajectory $\Phi(\tau,x)$ at $\tau\in[0,T]$ where
$d(\Phi(0,x),Y)<\varepsilon$
and $\varepsilon$ satisfies I-III. Suppose that
$x\in\partial N_{\tilde\delta}(\xi_{i_1})$ and the trajectory passes near
equilibria with the indices $[i_1,...,i_S]$. Denote
$d_s=d(g_{i_{s-1}}...g_{i_1}(x),Y)$ and $d_1=d(x,Y)$.
The proof employs the estimate
$d_{s+1}<AB(d_s)^{\rho_{i_s}}<(d_s)^{\tilde\rho_{i_s}}$, which follows
from III under the condition that $d_s<\tilde\varepsilon$.

For $d_2$ we have
$$d_2<AB(d_1)^{\rho_{i_1}}<(d_1)^{\tilde\rho_{i_1}}<
\tilde\varepsilon<\delta.$$
Similarly
$$d_{s+1}<AB(...AB(AB(d_1)^{\rho_{i_1}})^{\rho_{i_2}}...)^{\rho_{i_s}}<
(d_1)^{\tilde\rho_{i_s}...\tilde\rho_{i_1}}<\tilde\varepsilon<\delta,$$
as long as the indices $i_1,...,i_{s-1}$ are distinct.

If we have $i_s=i_d$ for some $1\le d\le s-2$, then the sequence of
equilibria $[\xi_{i_{d+1}},...,\xi_{i_s}]$ belongs to a heteroclinic cycle
$X\subset Y$. Hence
\begin{equation}\label{removeCycle}
d_{s+2}=d(g_{i_{s+1}}...g_{i_{d+1}}(g_{i_d}...g_{i_1}(x)),Y)<
(d_1)^{\tilde\rho(X)
\tilde\rho_{i_{s+1}}\tilde\rho_{i_d}...\tilde\rho_{i_1}}<
(d_1)^{\tilde\rho_{i_{s+1}}\tilde\rho_{i_d}...\tilde\rho_{i_1}}
<\tilde\varepsilon<\delta.
\end{equation}
We proceed as for the path for $i_{s+3}$, $i_{s+4}$, etc., provided that all indices
in the sequence  $i_1,...,i_d,i_{s+1},...,i_{s+k}$ are distinct. If
$i_{s+k}$ equals to some $i_p\in\{i_1,...,i_d,i_{s+1},...,i_{s+k-2}\}$, we
similarly remove the subsequence $\{i_{p+1},...,i_{s+k}\}$ and
obtain that
$$d_{s+k+2}<(d_1)^
{\tilde\rho_{i_1}...\tilde\rho_{i_p}\tilde\rho_{s+k+1}}
<\tilde\varepsilon<\delta.$$
Applying this procedure to the sequence of equilibria that is followed by
the trajectory $\Phi(\tau,x)$, we obtain that $d(\Phi(\tau,x),Y)<\delta$ for any
crossing of $N_{\tilde\delta}(\xi_j)$ at positive $\tau$.

To prove that a trajectory is attracted by $Y$, we recall that any trajectory
that is not attracted by $\xi_j\in Y$ passes near an infinite number of
equilibria.
Thus the sequence of equilibria visited by the trajectory must have repetitions.
The first repetition follows a cycle $X\subset Y$.
Due to our choice of $\varepsilon$ (condition II)
the estimate \eqref{removeCycle} may be refined to yield
$$d_{s+2}<(d_1)^{\tilde\rho(X)
\tilde\rho_{i_{s+1}}\tilde\rho_{i_d}...\tilde\rho_{i_1}}<
(d_1)^{(\tilde\rho(X)-1)
\tilde\rho_{i_{s+1}}\tilde\rho_{i_d}...\tilde\rho_{i_1}}
(d_1)^{\tilde\rho_{i_{s+1}}\tilde\rho_{i_d}...\tilde\rho_{i_1}}<
{1\over2}(d_1)^{\tilde\rho_{i_{s+1}}\tilde\rho_{i_d}...\tilde\rho_{i_1}}.$$
After removing of $k$ such subsequences we have
$$d_S<
{1\over2^k}(d_0)^{\tilde\rho_{i'}...\tilde\rho_{i''}},$$
where due to (III)
$(d_1)^{\tilde\rho_{i'}...\tilde\rho_{i''}}<\tilde\varepsilon$.
As $S\to\infty$ the number of removed cycles also tends to infinity we have $\lim_{S\to\infty}d_S=0$.
\qed

\begin{Cy}\label{CorMain}
Let $Y$ be an ac-network such that any heteroclinic cycle $X\subset Y$ satisfies
$\rho(X)>1$,
where the exponents $\rho_j$ are the ones derived in Lemmas
\ref{lem41}-\ref{lem45}. Then $Y$ is asymptotically stable.
\end{Cy}

Recalling that the quantities $\rho_j$ are defined in Lemmas \ref{lem41}, \ref{lem43} or
\ref{lem45}, the proof follows from the
property that $\rho_j>1$ only if $\xi_j$ is not an m-point and
Remark \ref{rem51}.

The conditions for asymptotic stability proven in Theorem \ref{th_as} are
not very tight,
in particular, in the employed estimates the exponents
$\rho_j$ for local maps do not take into account which faces a trajectory
enters and exits $N_{\tilde\delta}(\xi_j)$ through.
As can be seen from the proofs of Lemmas \ref{lem41}, \ref{lem43} and
\ref{lem45}, different values of $\rho_j$ can be obtained for different
prescribed entry and exit faces. Based on this, we propose   tighter
and more cumbersome, compared to Theorem \ref{th_as}, sufficient
conditions for asymptotic stability of a heteroclinic network.
We start by saying that a walk $\cP=[\xi_1\to...\to\xi_m]$ is
{\em semi-linear} if it does not have
repeating triplets $(\xi_{s-1}=\xi_{d-1},\xi_s=\xi_d,\xi_{s+1}=\xi_{d+1})$.

Decompose $\partial N_{\tilde\delta}(\xi_j)=\cup_{1\le s\le n}\cF_s$,
where $\cF_s$ is(are) the face(s) $|y_s|=\tilde\delta$ and $\by$ are the
local coordinates near $\xi_j$. Since we consider networks in $\R^n_+$,
there is just one face for each coordinate, except the radial one.
For a walk $\cP=[\xi_1\to...\to\xi_m]$ and an equilibrium $\xi_j\in\cP$
denote by $\cF_j^{in}$ the union of faces that intersect with the connection
$\xi_{j-1}\to\xi_j$ and by $\cF_j^{out}$ the faces that intersect with
the connection $\xi_j\to\xi_{j+1}$. For ac-networks
the sets $\cF_j^{in}$ and $\cF_j^{out}$ are unions of one or two faces,
depending on whether the respective connections are one- or two-dimensional.
Let $\Phi(\tau,x)$ enters $N_{\tilde\delta}(\xi_j)$
though $\cF_j^{in}$ and exits through $\cF_j^{out}$.
Then, by the same arguments that used to prove Lemmas  \ref{lem41}-\ref{lem45}, we have
the estimate
$$\td(\Phi(\tau^{out},x),Y)<A_j(\td(\Phi(\tau^{in},x),Y))^{\rho^*(\cP,\xi_j)},$$
with the exponent $\rho^*(\cP,\xi_j)\ge\rho_j$. The expressions for
$\rho^*(\cP,\xi_j)$ are similar to the ones for $\rho_j$ derived in
the lemmas and depend on dimensions of incoming and outgoing
connections and whether the eigenvectors are f-long/s-long vectors
of some $\Delta$-cliques. (To obtain $\rho^*(\cP,\xi_j)$ one just ignores
the parts of the lemmas that are not relevant to $\cF_j^{in}$ or $\cF_j^{out}$.)

\begin{theorem}\label{th_as2}
Let $Y$ be a clean heteroclinic network such that for any closed semi-linear
walk $\cP\subset Y$, $\cP=[\xi_1\to...\to\xi_m\to\xi_1]$, the local maps
$$
\phi(\cP,j):\,F_{\tilde\delta}(\xi_j,\kappa_{j-1,j})\to
F_{\tilde\delta}(\xi_j,\kappa_{j,j+1}),
$$
where $F_{\tilde\delta}(\xi_j,\kappa_{j-1,j})$ and
$F_{\tilde\delta}(\xi_j,\kappa_{j,j+1})$ are the faces of
$\partial N_{\tilde\delta}(\xi_j)$ that intersect with respective connections,
satisfy
$$
d(\Phi(\tau^{out},x),Y)<A_j(d(\Phi(\tau^{in},x),Y))^{\rho^*(\cP,j)}, \ A_j \in \R.
$$
Moreover, for any such walk
$$\rho^*(\cP)>1
\qquad\mbox{where}\qquad
\rho^*(\cP)=\prod_{1\le j\le m}\rho^*(\cP,j).
$$
Then $Y$ is asymptotically stable.
\end{theorem}

The proof is similar to the proof of Theorem \ref{th_as} and therefore is omitted.

\begin{remark}\label{stabtypes}
Any network discussed in part (i) of Theorem~\ref{th1} has a subcycle $X$,
such that any equilibrium $\xi_j\in X$ is an m-point for a $\Delta$-clique,
where the f-long and s-long vectors are the contracting and expanding
eigenvectors of $\xi_j$, respectively. This is also the
case for the networks discussed in Remark \ref{endcycl}. Instances of
such networks are given in Example~\ref{example5} in the next section.
Hence, by Lemma~\ref{lem42}, $\rho^*(X,\xi_j)=1$,
which implies that $\rho^*(X)=1$
and neither Theorem~\ref{th_as} nor Theorem~\ref{th_as2}
can be used to derive conditions for asymptotic stability of these networks.
Such networks were obtained in \cite{afraimovich16} as examples of
heteroclinic networks in the Lotka-Volterra system under the condition
that all equilibria have two-dimensional
unstable manifolds. As it was shown {\it ibid}
a sufficient condition of asymptotic stability of these networks is that
for each equilibrium $\xi_j\in Y$
$$|c_j/e_j|>1,\hbox{ where }c_j=\min_sc_{js}\hbox{ and }e_j=\min_qe_{jq},$$
$c_{js}$ are the contracting eigenvalues of $\rd f(\xi_j)$ and $e_{jq}$ are
the expanding ones.

For networks which were derived in part (ii) of Theorem~\ref{th1} and are
not the ones of Remark~\ref{endcycl}, there exists at least one $\xi_j$
such that the connection $[\xi_{j,1}\to\xi_{j+1,1}]$ does not exist.
Therefore, any semi-linear walk involves at least one of the connections
$[\xi_{j,s}\to\xi_{{j+1},0}]$ or $[\xi_{j,s}\to\xi_{j+2,0}]$, which is not
an f-long connection for any $\Delta$-clique.
Depending on the sequence of equilibria followed, either for $\xi_{{j+1},0}$ or
for $\xi_{{j+2},0}$, the exponent $\rho^*(\cP,\cdot)$ is calculated by
Lemma~\ref{lem41}. Therefore, by prescribing the respective expanding
eigenvalues to be sufficiently small, we can obtain $\rho^*(\cP)>1$
for any semi-linear walk in any such network. By any network we understang
any given grouping and set of connections allowed by part (ii) of Theorem~\ref{th1},
except for the ones of Remark~\ref{endcycl}.
\end{remark}

\section{Examples}\label{exa}

In this section we present five examples addressing the stability of six clean
heteroclinic networks, four of which are ac-networks.
The ac-networks appear in Examples~\ref{example1}, \ref{example2} and \ref{example5}.
For the networks in Examples~\ref{example1}, \ref{example3} and \ref{example4} sufficient conditions for asymptotic
stability are derived from Theorem~\ref{th_as}. The theorem, however, cannot
be employed to derive such conditions for the network in Example~\ref{example2}, because the
network contains a cycle $X$, where all equilibria $\xi_j\in X$ are m-points of
some $\Delta$-cliques. By Lemmas~\ref{lem43} and \ref{lem45} for
m-points $\rho_j\le1$, which implies that $\rho(X)\le1$. For this
network the sufficient conditions are derived from Theorem~\ref{th_as2}.
The last two networks are the ones discussed in the first part
of Remark~\ref{stabtypes}.

\begin{figure}[h]
\parbox{48mm}{
\begin{center}
\includegraphics[width=48mm]{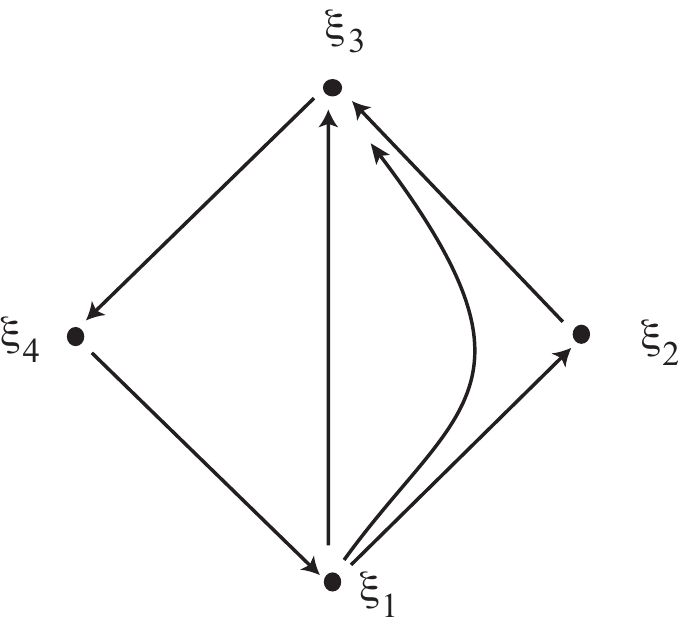}\\
(a)
\end{center}}
~~
\parbox{48mm}{
\begin{center}
\includegraphics[width=48mm]{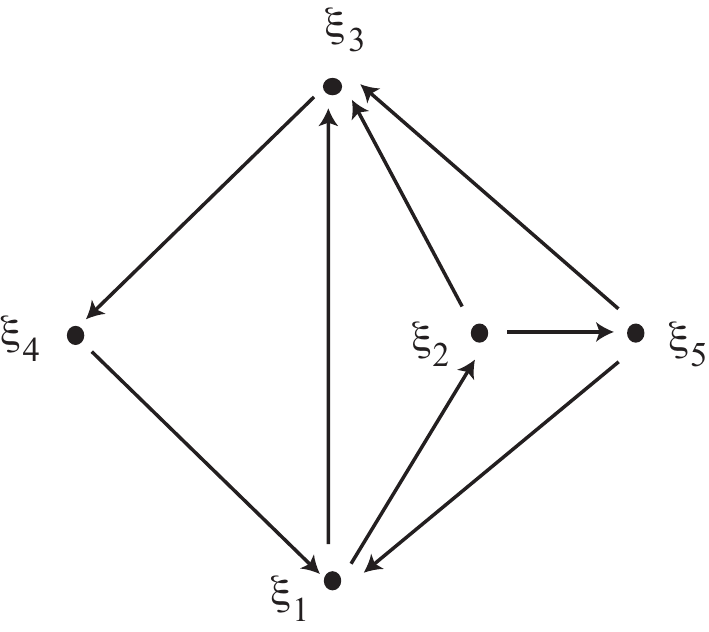}\\
(b)
\end{center}}
~~
\parbox{55mm}{
\begin{center}
\includegraphics[width=55mm]{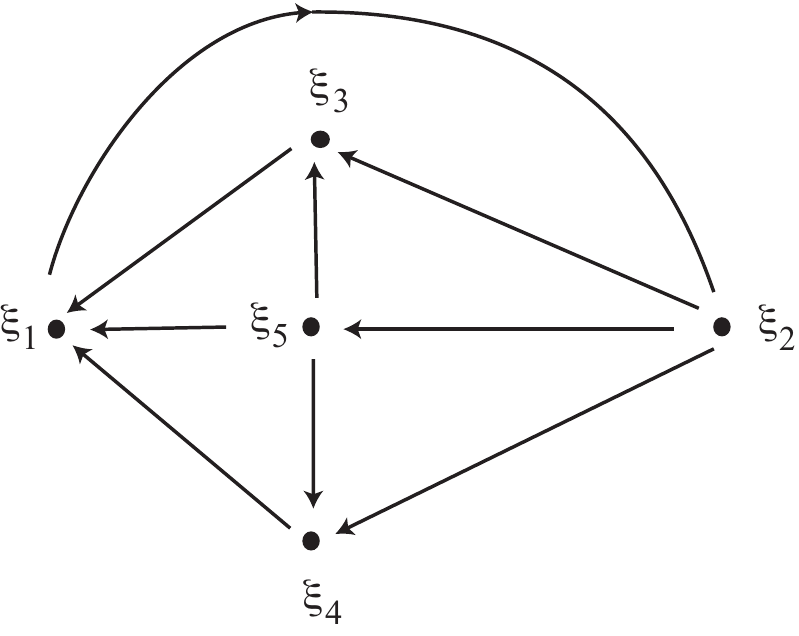}\\
(c)
\end{center}}
\caption{The networks of Examples~\ref{example1},\ref{example2} and \ref{example4}. \label{fig:Ex1}}
\end{figure}

\begin{example}\label{example1}
Consider a heteroclinic ac-network in $\R^4$ comprised of equilibria $\xi_j\in L_j$,
where $L_j$ are the coordinate axes, one-dimensional connections
$[\xi_3\to\xi_4]$
and $[\xi_4\to\xi_1]$ and a $\Delta$-clique  $\Delta_{123}$
(see Figure~\ref{fig:Ex1}a). By Lemmas \ref{lem41} and \ref{lem42} the values of $\rho_j$
for the equilibria are
$$
\rho_1=\displaystyle{{-\lambda_{14}\over\max(\lambda_{12},\lambda_{13})}},\quad
\rho_2=1,\quad
\rho_3=\displaystyle{{-\max(\lambda_{31},\lambda_{32})\over\lambda_{34}}},\quad
\rho_4=\displaystyle{\min\biggl({-\lambda_{43}\over\lambda_{41}},
1-{\lambda_{42}\over\lambda_{41}}\biggr)},
$$
where $\lambda_{ij}$ is the eigenvalue of $\rd f(\xi_i)$ associated with the
eigenvector $\be_j$.
The network is a union of two cycles, $[\xi_1\to\xi_2\to\xi_3\to\xi_4\to\xi_1]$
and $[\xi_1\to\xi_3\to\xi_4\to\xi_1]$, hence by Theorem \ref{th_as}
the sufficient condition for asymptotic stability of the network is
$${\lambda_{14}\max(\lambda_{31},\lambda_{32})
\over\lambda_{34}\min(\lambda_{12},\lambda_{13})}
\max\biggl({-\lambda_{43}\over\lambda_{41}},
1-{\lambda_{42}\over\lambda_{41}}\biggr)>1.$$
\end{example}

\begin{example}\label{example2}
Consider the ac-network $Y\subset\R^5$ shown in Figure~\ref{fig:Ex1}b. Its graph has
as a subset
the graph of the network considered in Example~\ref{example1}. It also has one more
equilirium $\xi_5\in L_5$ and three more connections: $[\xi_2\to\xi_5]$,
$[\xi_5\to\xi_1]$ and $[\xi_5\to\xi_3]$. Thus, the network
has three $\Delta$-cliques  $\Delta_{123}$, $\Delta_{253}$, and $\Delta_{513}$.
For the cycle
$X=[\xi_1\to\xi_2\to\xi_5\to\xi_1]$ the equilibria are m-points of some
$\Delta$-cliques, therefore $\rho(X)\le1$ by Lemmas \ref{lem43} and \ref{lem45}.
Hence, we cannot use Theorem \ref{th_as} to find conditions for
asymptotic stability of $Y$.
The values of $\rho^*(\cP)$ calculated for all semi-linear paths
$\cP\subset Y$ are presented in the table below:
\begin{center}
$$
\begin{array}{l|l}
\cP & \rho^*(\cP)\\
\hline
(\xi_1\to\xi_2\to\xi_5\to\xi_1)& q_{125}q_{251}q_{512}\\
(\xi_4\to\xi_1\to\xi_3\to\xi_4)& q_{413}q_{134}q_{341}\\
(\xi_4\to\xi_1\to\xi_2\to\xi_3\to\xi_4)& q_{412}q_{123}q_{234}q_{341}\\
(\xi_4\to\xi_1\to\xi_2\to\xi_5\to\xi_3\to\xi_4)&
q_{412}q_{125}q_{253}q_{534}q_{341}\\
(\xi_4\to\xi_1\to\xi_2\to\xi_5\to\xi_1\to\xi_3\to\xi_4)&
q_{412}q_{125}q_{251}q_{513}q_{134}q_{341}\\
(\xi_4\to\xi_1\to\xi_2\to\xi_5\to\xi_1\to\xi_2\to\xi_3\to\xi_4)&
q_{412}q_{125}q_{251}q_{512}q_{123}q_{234}q_{341}\\
(\xi_4\to\xi_1\to\xi_2\to\xi_5\to\xi_1\to\xi_2\to\xi_5\to\xi_3\to\xi_4)&
q_{412}q_{125}q_{251}q_{512}q_{125}q_{253}q_{534}q_{341}\\
(\xi_4\to\xi_1\to\xi_2\to\xi_5\to\xi_1\to\xi_2\to\xi_5\to\xi_1\to\xi_3\to\xi_4)&
q_{412}q_{125}q_{251}q_{512}q_{125}q_{251}q_{513}q_{134}q_{341}
\end{array}
$$
\end{center}
where $q_{ijk}$ denotes $\rho^*(\cP,\xi_j)$, i.e., the exponent of the local
map near $\xi_j$ with the contracting
eigenvector $\bv_i$ and the expanding $\bv_k$. They were calculated
from Lemmas \ref{lem41}-\ref{lem42}:
$$
\begin{array}{l}
q_{412}=\displaystyle{\min\biggl({-\lambda_{14}\over\lambda_{12}},
1-{\lambda_{15}\over\lambda_{12}}\biggr)},\quad
q_{413}=\displaystyle{\min\biggl({-\lambda_{14}\over\lambda_{13}},
1-{\lambda_{15}\over\lambda_{13}}\biggr)},\quad
q_{512}=\displaystyle{\min\biggl({-\lambda_{15}\over\lambda_{12}},
1-{\lambda_{14}\over\lambda_{12}}\biggr)},\\
q_{513}=\displaystyle{{\lambda_{15}\over\lambda_{15}-\lambda_{13}}},\quad
q_{123}=\displaystyle{{\lambda_{21}\over\lambda_{21}-\lambda_{23}}},\quad
q_{125}=\displaystyle{\min\biggl({-\lambda_{21}\over\lambda_{25}},
1-{\lambda_{24}\over\lambda_{25}}\biggr)},\\
q_{134}=\displaystyle{\min\biggl({-\lambda_{31}\over\lambda_{34}},
1-{\max(\lambda_{32},\lambda_{35})\over\lambda_{34}}\biggr)},\quad
q_{234}=\displaystyle{\min\biggl({-\lambda_{32}\over\lambda_{34}},
1-{\max(\lambda_{31},\lambda_{35})\over\lambda_{34}}\biggr)}\\
q_{534}=\displaystyle{\min\biggl({-\lambda_{35}\over\lambda_{34}},
1-{\max(\lambda_{31},\lambda_{32})\over\lambda_{34}}\biggr)},\quad
q_{341}=\displaystyle{\min\biggl({-\lambda_{43}\over\lambda_{41}},
1-{\max(\lambda_{42},\lambda_{45})\over\lambda_{41}}\biggr)},\\
q_{251}=\displaystyle{\min\biggl({-\lambda_{52}\over\lambda_{51}},
1-{\lambda_{54}\over\lambda_{51}}\biggr)},\quad
q_{253}=\displaystyle{{\lambda_{52}\over\lambda_{52}-\lambda_{51}}}.
\end{array}
$$
By Theorem \ref{th_as2} a sufficient condition for the stability of the network is
that $\rho^*(\cP)>1$ for all $\cP$ listed in the table. (These inequalities
hold true, e.g., if $|\lambda_{42}|$, $|\lambda_{43}|$ and $|\lambda_{45}|$ are significantly larger than other
eigenvalues and in addition for the cycle $[\xi_1\to\xi_2\to\xi_5\to\xi_1]$ the
contracting eigenvalues are larger than the expanding ones, in absolute value.)
\end{example}

\begin{example}\label{example3}
The network shown in Figure~\ref{fig2}a was discussed in \cite{AshCasLoh18}
(Subsection 4.1
{\it ibid}) as an example of a clean heteroclinic network in $\R^4$.
The equilibria $\xi_1,\xi_2,\xi_3,\xi_4$ belong to the respective coordinate
axes, while $\xi_5$ belongs to the coordinate plane $P_{34}$.
This is not an ac-network because $\xi_5$ does not belong to a coordinate axis.
However, since it is clean, we can apply Theorem \ref{th_as}
to derive conditions for asymptotic stability of this network. The values of
$\rho_j$ are:
$$
\begin{array}{l}
\rho_1=\displaystyle{{-\max(\lambda_{13},\lambda_{14})\over\lambda_{12}}},\quad
\rho_2=\displaystyle{{-\lambda_{21}\over\max(\lambda_{23},\lambda_{24})}},\quad
\rho_3=\displaystyle{{-\max(\lambda_{32},\lambda_{34})\over\lambda_{31}}},\\
\rho_4=\displaystyle{{-\max(\lambda_{42},\lambda_{43})\over\lambda_{41}}},\quad
\rho_5=1.
\end{array}
$$
Applying the theorem to the four cycles the network is comprised of, we obtain
that the conditions for asymptotic stability are:
$$\rho_1\rho_2\min(\rho_3,\rho_4)>1.$$
The network in Figure~\ref{fig2}b can be obtained from the one in
Figure~\ref{fig2}a by removing $\xi_5$ and the
connections $[\xi_2\to\xi_5]$, $[\xi_5\to\xi_3]$ and $[\xi_5\to\xi_4]$. Under the condition
that the constants of the global maps near the network are uniformly bounded,
the above condition is also applicable to this network.
\end{example}

\begin{example}\label{example4}
The network shown in Figure~\ref{fig:Ex1}c,
the stability and bifurcations which were studied in \cite{kpr},
can be obtained from the one in
Figure~\ref{fig2}a by adding a two-dimensional connection $[\xi_5\to\xi_1]$.
Since not all equilibria are on coordinate axes this not an ac-network.
Below we apply apply Theorem \ref{th_as} to this network.
The values of $\rho_j$ are:
$$
\begin{array}{l}
\rho_1=\displaystyle{{-\max(\lambda_{13},\lambda_{14})\over\lambda_{12}}},\quad
\rho_2=\displaystyle{{-\lambda_{21}\over\max(\lambda_{23},\lambda_{24})}},\quad
\rho_3=\displaystyle{\min\biggl({-\lambda_{32}\over\lambda_{31}},1\biggr)},\\
\rho_3=\displaystyle{\min\biggl({-\lambda_{42}\over\lambda_{41}},1\biggr)},\quad
\rho_5=\displaystyle{{\lambda_{52}\over\lambda_{52}-\lambda_{51}}}.
\end{array}
$$
The network has five subcycles, which implies that by Theorem~\ref{th_as}
the conditions for stability are:
$$\rho_1\rho_2\rho_3>1,\quad \rho_1\rho_2\rho_4>1,\quad \rho_1\rho_2\rho_5>1,
\quad \rho_1\rho_2\rho_5\rho_3>1,\quad \rho_1\rho_2\rho_5\rho_4>1,$$
which can be combined into
$$\rho_1\rho_2\rho_5\min(\rho_3,\rho_4,1)>1.$$
\end{example}

\begin{figure}[h]
\vspace*{5mm}
\hspace*{25mm}\includegraphics[width=55mm]{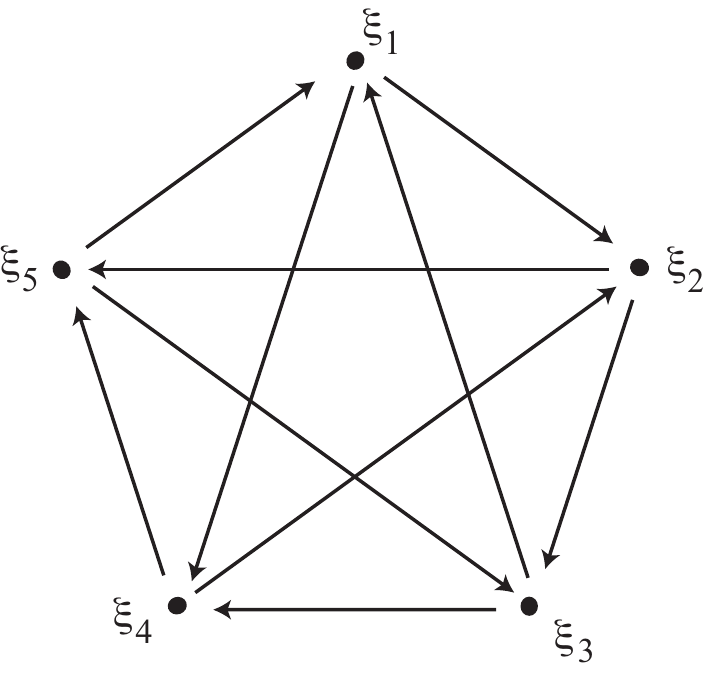}~~
\hspace*{15mm}\includegraphics[width=62mm]{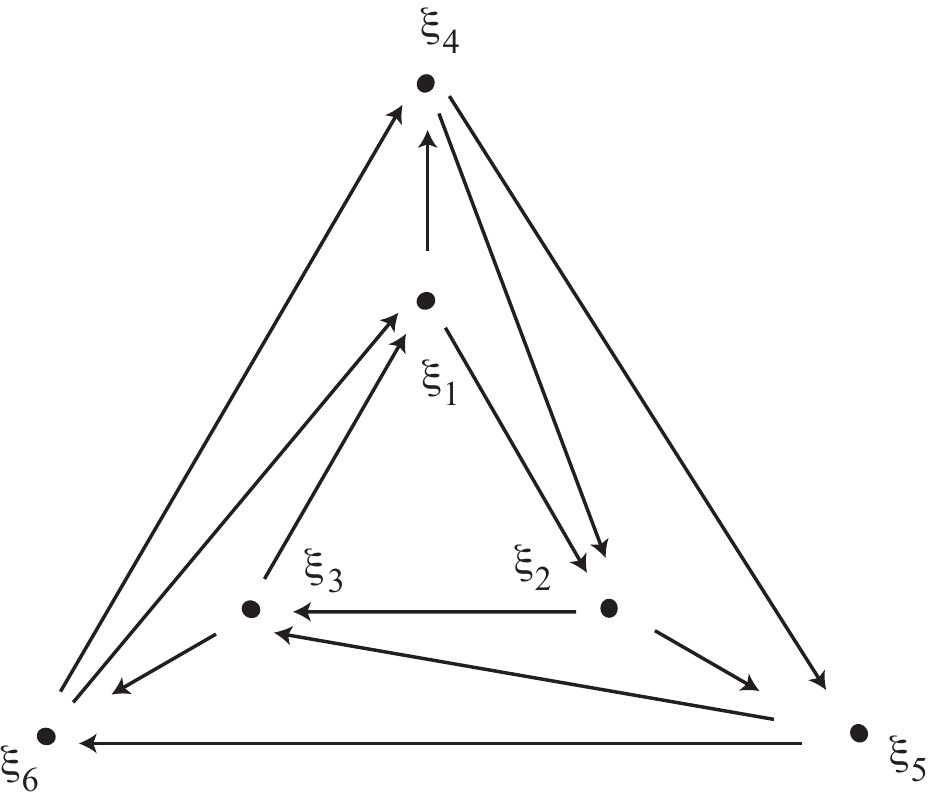}

\caption{The networks of Example~\ref{example5}: $Y_5$ on the left and $Y_6$ on the right. \label{fig:Ex5}}
\end{figure}

\begin{example}\label{example5}
The ac-networks $Y_5\subset\R^5$ and $Y_6\subset\R^6$ shown in Figure~\ref{fig:Ex5}
are the ones discussed in Remark~\ref{stabtypes}.
They have subcycles where all equilibria have f-long vectors
as contracting eigenvectors and s-long vectors of the same $\Delta$-clique
as expanding eigenvectors. The subcycles are
$$
Y_5\supset X_5=[\xi_1\to\xi_4\to\xi_2\to\xi_5\to\xi_3\to\xi_1]
$$
and
$$
Y_6\supset X_6=[\xi_1\to\xi_4\to\xi_2\to\xi_5\to\xi_3\to\xi_6\to\xi_1],
$$
which implies that $\rho^*(X_5)=\rho^*(X_6)=1$.
Equality prevents us from using  Theorems~\ref{th_as} and  \ref{th_as2}.
As we noted, according to
\cite{afraimovich16}, for such cycles in the Lotka-Volterra system the
sufficient conditions for asymptotic stability are that $|c_j/e_j|>1$
for any $\xi_j\in Y$.
\end{example}

\section{Concluding remarks}\label{conc}

In this paper we introduce ac-networks and prove sufficient conditions for their
asymptotic stability, which are applicable
also in more
general cases. Using these results, one can construct various kinds of
dynamical systems that have asymptotically stable heteroclinic networks,
for example, using
a generalisation of the methods proposed in
\cite{AshPos13} or \cite{clp18,pl19}. Alternatively,
given any graph of the structure described in Theorem \ref{th1}, there
exists a Lotka-Volterra system that has an associated network. Conditions
on the coefficients of the system that ensure the existence of
heteroclinic trajectories can be found, e.g., in \cite{op10}. Note also that an
ac-network can be realised as a network on a simplex or in coupled
identical cell system \cite{fie17}.

The conditions for asymptotic stability that we derive are general and
have a simple form. This is related to the fact that they are not accurate.
In \cite{pcl18} we obtained necessary and sufficient conditions for
fragmentary asymptotic stability of a particular heteroclinic network in
$\R^6$, which involve rather cumbersome expressions. The study indicated
that derivation of conditions for stability of a network in $\R^n$ is a
highly non-trivial task, unless in some special cases or if $n$ small.
Here simplicity is achieved by compromising in accuracy. It would be of
interest to compare the sufficient conditions that we obtained with
precise conditions for asymptotic stability that could be derived for some
ac-networks.

The present study is also a first step towards considering more general types
of heteroclinic networks, such as networks involving unstable manifolds of
dimension greater than two, with equilibria not just on coordinate axes, with
more than one equilibrium per coordinate axis, networks that are not clean,
as well as networks in systems equivariant under the action of symmetry groups
other than $\Z^n$.

For the classification of the
ac-networks (Theorem~\ref{th1}) the requirements that the equilibria involved in the network belong to coordinate
axes and $\Z^n$-equivariance are essential.
This is not the case of the lemmas and theorems
in Section~\ref{sufcond} where the conditions are more general. Note that these
results are applied to networks considered in Examples~\ref{example3} and \ref{example4}
 in Section~\ref{exa} that are not ac-networks.
Our results on asymptotic stability are applicable to any robust heteroclinic
network where (i) the unstable manifolds of equilibria are one- or
two-dimensional; (ii) the global maps between neighbourhoods are bounded (e.g.
due to the compactness of the network); (iii) the robust distance can be used
instead of the distance based on the maximum (or any other equivalent) norm.
Apparently,  when
the unstable manifold of an equilibrium has dimension three or higher an
estimate for the local maps similar to the ones derived in Subsection~\ref{locglob}
holds true with the exponent depending on the eigenvalues of the Jacobian.
However, we expect the dependence to be more complex than the one obtained
for two-dimensional unstable manifolds.

\section*{Acknowledgements}
O.P. had partial financial support  by a grant from project STRIDE [NORTE-01-0145-FEDER-000033] funded by FEDER --- NORTE 2020 --- Portugal.
 S.C. and I.L. were partially supported by CMUP (UID/MAT/00144/2019), which is funded by FCT (Portugal) with national (MCTES) and European structural funds through the programs FEDER, under the partnership agreement PT2020.

\end{document}